\documentclass[centertags,leqno,11pt]{article}

\usepackage{amsmath}
\usepackage{amssymb}
\usepackage{latexsym}
\usepackage{amsxtra}
\usepackage{amscd}
\usepackage{theorem}

\usepackage{graphics}
\usepackage{epic}

\theoremstyle{change}

{\theorembodyfont{\slshape}

\newtheorem{thm}{Theorem.}[section]
\newtheorem{cor}[thm]{Corollary.}
\newtheorem{lem}[thm]{Lemma.}
\newtheorem{prop}[thm]{Proposition.}

\newtheorem{defn}[thm]{Definition.}}

{\theorembodyfont{\rmfamily}

\newtheorem{rem}[thm]{Remark.}

}

\renewcommand{\em}{\sl}

\newcommand{\proof}{\noindent {\bf Proof:\ }}
\newcommand{\Endproof}{\hspace*{\fill} $\Box$ \vspace{1ex} \noindent }

\makeatletter
\renewcommand{\subsection}{\@startsection{subsection}{2}%
{\z@}{-3.25ex plus -1ex minus-.2ex}{-1em}{\bf}} \makeatother

\newcommand{\PP}{\mathbb{P}}
\newcommand{\ZZ}{\mathbb{Z}}
\newcommand{\CC}{\mathbb{C}}
\newcommand{\RR}{\mathbb{R}}
\newcommand{\QQ}{\mathbb{Q}}
\newcommand{\NN}{\mathbb{N}}
\newcommand{\FF}{\mathbb{F}}
\renewcommand{\AA}{\mathbb{A}}
\newcommand{\GG}{\mathbb{G}}
\newcommand{\F}{\mathcal{F}}

\newcommand{\V}{\mathcal{V}}

\newcommand{\K}{\mathcal{K}}

\newcommand{\OO}{\mathcal{O}}

\newcommand{\eps}{{\epsilon}}

\newcommand{\geom}{{\rm geo}}

\newcommand{\bQl}{{\bar{\QQ}_\ell}}
\newcommand{\naive}{{\rm naive}}
\newcommand{\Sp}{{\rm Sp}}

\newcommand{\Trace}{{\rm Tr}}

\newcommand{\g}{{\bf g}}

\newcommand{\GL}{{\rm GL}}
\newcommand{\SL}{{\rm SL}}

\newcommand{\PSL}{{\rm PSL}}
\newcommand{\Gal}{{\rm Gal}}
\newcommand{\Aut}{{\rm Aut}}
\newcommand{\Spec}{{\rm Spec\,}}

\newcommand{\T}{{\cal T}}
\renewcommand{\middle}{{\rm mid}}

\newcommand{\MC}{{\rm MC}}

\newcommand{\an}{{\rm an}}

\newcommand{\et}{{\rm\acute{e}t}}

\newcommand{\LS}{{\rm LS}}

\renewcommand{\P}{{\cal P}}

\newcommand{\To}{\;\longrightarrow\;}

\newcommand{\Mapsto}{\;\longmapsto\;}

\newcommand{\diag}{{\rm diag}}

\newcommand{\im}{{\rm Im}}

\newcommand{\al}{{\alpha}}

\newcommand{\D}{{\rm D}}
\newcommand{\uo}{{{\bf u}}}

\newcommand{\HHH}{{{\cal H}}}

\newcommand{\vo}{{{\bf v}}}

\newcommand{\q}{{{\rm d}}}

\newcommand{\pr}{{{\rm pr}}}
\newcommand{\opr}{{\overline{\rm pr}}}
\newcommand{\tim}{{\boxtimes}}

\newcommand{\rk}{{\rm rk}}

\newcommand{\J}{{\rm J}}

\newcommand{\Constr}{{\rm Constr}}

\newcommand{\Frob}{{\rm Frob}}

\newcommand{\chara}{{\rm char}}

\newcommand{\Perv}{{\rm Perv}}

\newcommand{\FFF}{{\cal F}}\newcommand{\GGG}{{\cal G}}

\newcommand{\End}{{\rm  End}}

\newcommand{\SO}{{\rm SO}}
\renewcommand{\char}{{\rm char}}
\renewcommand{\H}{{\cal H}}
\renewcommand{\L}{{\cal L}}

\newcommand{\tame}{{\rm tm}}
\renewcommand{\d}{{\rm d}}

\numberwithin{equation}{subsection}
\numberwithin{thm}{subsection}
\theoremstyle{plain}

\title{On the middle convolution of local systems. \\
{\large With an Appendix by
Stefan Reiter and Michael Dettweiler}.}

\author{
   Michael Dettweiler \\ IWR Heidelberg
  }

\begin{document}
\maketitle 
\begin{abstract} 
We study the  middle convolution
of local systems in the setting of singular and \'etale 
cohomology.  
We give a motivic
interpretation of the middle convolution
in the \'etale case and prove an independence-of-$\ell$-result  which
 yields a description of 
the determinant. 
We employ these methods to realize  special linear groups regularly as
Galois groups over $\QQ(t).$ In an appendix to this article, 
written jointly with S. Reiter, we prove
the existence of a new motivic local system whose monodromy is
dense in the exceptional simple group of type~$G_2.$ 
\end{abstract}

\tableofcontents 

%------------------------------------------------------

%\addcontentsline{toc}{section}{Introduction}
%\setcounter{page}{1}
\section*{Introduction}\label{Introduction}

On the affine line $\AA^1$ over an algebraically closed 
field or a finite field $k$ one has 
 the 
derived category $\D^b_c(\AA^1,\bar{\QQ}_\ell)$ of 
constructible $\bar{\QQ}_\ell$-sheaves with 
bounded cohomology, cf. \cite{BBD} and \cite{KiehlWeissauer}. 
This category supports Grothendieck's 
six operations, so one has the notion of 
a higher direct image and a higher direct image 
with compact support.  
Consider the addition map 
$$ \pi: \AA^1\times \AA^1 \To \AA^1.$$ 
If $\K$ and $\L$ are elements in $\D^b_c(\AA^1,\bar{\QQ}_\ell),$
  then one can consider 
their {\it $*$-convolution} as the higher direct image 
$$ \K\ast_* \L:=R\pi_*(\K \boxtimes \L)\in \D^b_c(\AA^1,\bar{\QQ}_\ell),$$
where $\K\boxtimes \L\in \D^b_c(\AA^1\times \AA^1,\bar{\QQ}_\ell)$ denotes the external tensor 
product of $\K$ and $\L.$ One can also define the $!$-convolution
of $\K$ and $\L$ as the higher direct image with compact supports
$$ \K\ast_! \L:=R\pi_!(\K \boxtimes \L)\in \D^b_c(\AA^1,\bar{\QQ}_\ell).$$
For more details, the reader should consult
 the book {\it Rigid Local Systems}
by N.~Katz~\cite{Katz96}.
The structure of the convolution is very complicated in general,
so it is convenient to restrict the above construction 
to smaller subcategories of $\D^b_c(\AA^1,\bar{\QQ}_\ell).$ 
A natural candidate to work with is 
the abelian category of perverse sheaves 
$\Perv(\AA^1,\bar{\QQ}_\ell) \subseteq \D^b_c(\AA^1,\bar{\QQ}_\ell),$ cf.~\cite{BBD} and 
\cite{KiehlWeissauer}.  
This category is not stable under convolution
in general, but we can consider 
 perverse sheaves  $\K\in \Perv(\AA^1,\bar{\QQ}_\ell)$ which
have  
the property that for any other $\L\in \Perv(\AA^1,\bar{\QQ}_\ell),$ the convolutions
$\K\ast_! \L$ and $ \K\ast_* \L$ are perverse. In this case 
we say that $\K$ has the property $\P.$ 
For example, an object $\K\in \Perv(\AA^1,\bar{\QQ}_\ell)$ has the 
property $\P$ if it is irreducible 
and if its  isomorphism class 
is not translation invariant, see 
\cite{Katz96}, Cor.~2.6.10. 

If $\K\in \Perv(\AA^1,\bar{\QQ}_\ell)$ has the 
property $\P,$ then  one can define the 
{\it middle convolution} of $\K$ and $\L\in \Perv(\AA^1,\bar{\QQ}_\ell)$ 
as 
\begin{equation}\label{eqimmid}
 \K*_\middle\L =\K*\L=\im\left(\K\ast_! \L \to \K\ast_* \L\right),\end{equation}
see \cite{Katz96}, Chap.~2.6 (we mostly omit the subscript 
${}_\middle$ for notational reasons). This is a perverse sheaf which 
has again 
the property $\P,$ 
see \cite{Katz96}, Cor.~2.6.17. Let
 $$j:\AA^1_x\times \AA^1_y \hookrightarrow \PP^1_x\times \AA^1_y,\quad 
\opr_2:\PP^1_x\times \AA^1_y\to \AA^1_y,$$ 
and 
$$ \d:\AA^1_x\times \AA^1_y\to \AA^1_{y-x},\, (x,y)\mapsto y-x.$$
It is more or less tautological that 
 the middle convolution 
can then be interpreted in terms of middle direct images 
 as follows (see \cite{Katz96}, Prop.~2.8.4):
\begin{equation}\label{eqimmid2}
 \K*\L =R \,\opr_{2*} \left(j_{!*}(\pr_1^*\K\otimes \d^*\L)\right),\end{equation}
where $j_{!*}(\pr_1^*\K\otimes \d^*\L)$ denotes the middle direct image
(sometimes called intermediate extension) of $\pr_1^*\K\otimes \d^*\L.$
One reason why one is interested in the middle 
 convolution is that 
 $\K\ast \L$ is 
often irreducible, while the convolutions
$\K\ast_* \L$ and $\K\ast_! \L$
 are not irreducible. 
A striking application of  the concept of middle 
convolution is Katz' existence algorithm for irreducible 
rigid local systems, see 
\cite{Katz96}, Chap.~6.\\

In analogy to Formula \eqref{eqimmid2},
one can give a definition of 
the middle convolution of local systems and of 
lisse \'etale sheaves
(sometimes called \'etale local systems), simply by 
working with suitably small open subsets of 
$\AA^2$ and $\AA^1$ over which all 
constructions yield again local systems, see Sections~\ref{sectconvdef} and 
\ref{secvondefetale}.  Our approach is strongly motivated
by the definitions and the results of 
\cite{Katz96}, Chap.~8. 

The restriction of the 
middle convolution of perverse sheaves to a suitably small 
dense open subset of $\AA^1$ is often given by the middle convolution of 
local systems, see Section~\ref{secperverse}. Therefore,
the concept of middle convolution of local systems
can also be seen as a 
powerful tool for the study of the monodromy  of middle
convoluted perverse sheaves.
 
We then give some applications to the
 inverse Galois problem and  the construction problem 
of interesting motivic local systems.  
The philosophy is that by convoluting elementary objects, 
like the local systems  associated to cyclic or 
dihedral Galois covers of open subsets $S\subseteq \AA^1,$  one obtains highly non-trivial
  Galois representations and motives.\\

Let us now describe the content of this paper in more detail: 
We start in Section~\ref{sectconvdef1} by  working on the 
analytic affine line  
$\AA^1(\CC).$ 
In Section~\ref{sectconvdef} we give the definition of 
the middle convolution of local systems 
as the parabolic cohomology 
of a variation of local systems. 
 The main results  are the following: 
\begin{itemize}
\item We  derive a formula for
the rank
of the middle convolution
$\V_1\ast \V_2,$
see Prop.~\ref{dimensione}. 
\item We study the important case of convoluting local systems with Kummer 
sheaves and relate the monodromy of such convolutions to the
tuple transformation $\MC_\lambda$ of \cite{dr00} and \cite{dr03}. 
\item An 
irreducibility criterion for the convolution of 
some local systems is given in  Thm.~\ref{thmirrd}. 
\item The effect of the middle convolution on 
the local monodromy is determined in some cases, 
see Section~\ref{secmonodrconv}. 
\end{itemize}

In Section~\ref{secetaleconN}, 
we  study the middle convolution of lisse \'etale sheaves.
We derive the following results:
\begin{itemize}
\item If the ground ring $R$ is contained in $\CC,$ then 
the geometric monodromy of the middle convolution of
lisse sheaves can be computed using 
the concept of analytification, see Prop.~\ref{dimension}. 
\item We show how the middle convolution of perverse sheaves
 is related to the middle convolution of lisse sheaves, see 
Section~\ref{secperverse}.
\item 
Analogous to \cite{Katz96}, Thm.~5.5.4, we prove 
that  independence-of-$\ell$ is preserved by convolution, see Thm.~\ref{indep}.
The theory of Hecke characters and a result of Henniart~\cite{Henniart81} imply  then that over $\QQ,$ the determinant of the middle convolution
 is the product of the geometric determinant with 
 a finite character and a power of the cyclotomic character.
\item We give a motivic interpretation of 
the middle convolution, see Thm.~\ref{thmpropkat}.  This is close
to the results in \cite{Katz96}, Chap.~8 (but not contained in them). \\
\end{itemize}

We then give two applications of the above methods.
The first one relies on the irreducibility criterion of Thm.~\ref{thmirrd}
and the above mentioned application of the theory of Hecke characters
 to the determinant (see Cor.~\ref{corsldreia}): \\

\noindent {\bf Theorem~II.} {\it Let $\FF_q$ be the finite field 
of order $q=\ell^k,$ where $k\in \NN .$ 
Then the special linear group $\SL_{2n+1}(\FF_q)$  
occurs regularly as Galois
group over $\QQ(t)$ if $$q\equiv 5\mod 8\quad \mbox{and}\quad 
n> 6+2\varphi((q-1)/4)$$ ($\varphi$ denoting
Euler's $\varphi$-function){\it .}}\\

The theorem implies that, under the above restrictions,
the simple  groups 
$\PSL_{2n+1}(\FF_q)$ occur regularly as Galois
group over $\QQ(t).$ 
The latter result is the first result on regular Galois
realizations of the groups ${\rm PSL}_n(\FF_q)$ over $\QQ(t),$ 
 where 
$$(n,q-1)=[{\rm PGL}_n(\FF_q):{\rm PSL}_n(\FF_q)]> 2\,.$$
We also realize the underlying profinite special linear groups 
regularly as Galois groups over $\QQ(t)$ (Thm.~\ref{thmrealierung1}). \\

In the  appendix to this article, written 
jointly by S.~Reiter and the author,  we prove the existence 
of a new motivic lisse sheaf $\H$ whose monodromy is dense in
the exceptional algebraic group of type $G_2.$ This result relies
on our results on 
the monodromy of the middle convolution with Kummer 
sheaves
(Section~\ref{secirrkumm}) and on
the motivic interpretation  
of the middle convolution (Section~\ref{secmotiv}).
 The lisse sheaf $\H$ is 
not rigid, contrary to the local systems considered
in \cite{DK}, \cite{DR07}. The explicit determination
of the monodromy seems to be necessary in this case, since there exist
other local systems with the same local monodromy whose 
monodromy is dense in the special orthogonal group $\SO_7$ 
(see the remark at the end of the appendix). 
It seems to be remarkable, that the 
weight of (a certain arithmetic extension of) $\H$  is $4,$ contrary
to the systems of \cite{DR07}, \cite{DK}, where the 
weight is $6.$\\

{\bf Acknowledgments:} This article presents a revised part of my 
Habilitation Thesis (Heidelberg, 2005). 
I  heartily thank  B.H. Matzat for 
his support during the last years
 and  for many valuable discussions 
on the subject of this work. 
I thank M. Berkenbosch, D. Haran
J. Hartmann, M. Jarden, U. K\"uhn,  S. Reiter, A. R\"oscheisen and S. Wewers 
for valuable comments and discussions. 
Part of this work was written 
during my stays at the School of Mathematics of the Tel Aviv University
(spring 2004), the Institute for Advanced Study (IAS)
in Princetion (spring 2007), and at the Insitut des Hautes \'Etudes
Scientifiques (IHES), Bures sur Yvette, Paris (Fall 2008). I thank these institutes 
for creating a very friendly  and inspiring atmosphere. 

The appendix has grown out from a question of P. Deligne about 
a possible motivic interpretation of the 
other $G_2$-rigid (but non-$\GL_7$-rigid) local systems that exist 
besides the $G_2$-rigid local systems 
 considered in \cite{DR07}. The authors are indepted to Professors P. Deligne 
and N. Katz for their interest and several valuable remarks and discussions
on the subject.

\section{Convolution of local systems}\label{sectconvdef1}

Throughout this section we
will write $\AA^1,\, \PP^1 , \ldots$ instead of $\AA^1(\CC),\,
\PP^1(\CC), \ldots$ and  view these objects equipped with their
associated topological and complex analytic structures.
 Let $X$ be a connected topological manifold.
 The multiplication in the fundamental group $\pi_1(X,x)$
is induced by the path product for which we use 
the following convention: Let  $\gamma, \gamma'$ be two 
closed paths at $x\in X.$ Then their product $\gamma \gamma'$ is given
by first walking along $\gamma'$ and then walking along $\gamma$. 
Endomorphisms of vector spaces act from the right throughout the
article. 

 \subsection{Some notation.}\label{braidgroups}
Let $U_0= \AA^1\setminus \uo,$
where $\uo:= \{ u_1,\ldots,u_r\}$ is a finite subset 
of $ \AA^1.$   It is well known that
there exist generators $\al_1,\ldots,\al_{r+1}$ of $\pi_1(U_0,u_0)$
which are the homotopy classes of simple loops 
around the points $u_1,\ldots,u_r,\infty$ (resp.) and 
which satisfy the product relation $\al_1\cdots\al_{r+1}=1.$ 
Let $R$ be a commutative ring with a unit and let 
$X$ be a connected topological manifold. A {\em local system of $R$-modules}
is a sheaf $\V$ which is locally isomorphic to ${R}^n$ for
some $n\in \NN.$  The
number $n$ is called the {\em rank} of $\V$ and is denoted by
$\rk(\V).$ Let $\LS_R(X)$ denote the category of local systems of
$R$-modules on $X.$ Any local system $\V \in \LS_R(X)$ gives rise
to its {\em monodromy representation}
$$ \rho_\V:\pi_1(X,x) \To \GL(V), \, \gamma \Mapsto \rho_\V(\gamma),$$
where $V =\V_{x}$ denotes the stalk of $\V$ at $x,$ 
cf.~\cite{Deligne70}.

\begin{defn}\label{deff}
{\rm Let 
$U_0:=\AA^1\setminus \uo$ be as above 
 and fix generators $\al_1,\ldots,\al_{r+1}$ of
$\pi_1(U_0,u_0)$ as above.
Let $\V$ be a local system on $U_0$ with 
 monodromy representation $\rho_\V.$ Then $\rho_\V$ is 
uniquely determined by
 the {\em monodromy tuple} of $\V$ (with respect to 
$\alpha_1,\ldots\alpha_r$):
$$T_\V=(T_1,\ldots,T_{r+1})\in \GL(V)^{r+1},\quad 
 T_i=\rho_\V(\alpha_i).
$$}
\end{defn}
The following definition is motivated by the results in \cite{Katz96}, Chap.~5:
\begin{defn}\label{defconvosheaf}{\rm  
A local system $\V \in \LS_K(U_0)$ has the property $T,$  
if it is irreducible and if there exist 
at least two components $T_i,\,T_j, \, i<j\leq r,$ 
of the monodromy tuple $\T_\V=(T_1,\ldots
T_r,T_{r+1})\in \GL(V)$ which are not the identity. Let $\T_K(U_0)$ be the category
of local systems in $\LS_K(U_0)$ having the property $T.$}
\end{defn}

\begin{defn} {\rm If $\V$ is a local system on $U_0$ and if $j:U_0\to \PP^1$ is the 
natural inclusion, then the {\it parabolic cohomology group }
of $\V$ is defined to be
$ H^1_p(U_0,\V):=H^1(\PP^1,j_*\V),$ cf.~\cite{dw03}, Section~1. }
\end{defn}

\subsection{The middle  convolution
of local systems.}\label{sectconvdef}

For $\uo :=\{ x_1,\ldots,x_p\}\subseteq \AA^1$ and $\vo
:=\{y_1,\ldots,y_q\}\subseteq \AA^1,$ set
$$ \uo \ast \vo := \{ x_i + y_j \mid i=1,\ldots,p, \,\, j = 1,\ldots
,q \}\,\subseteq \AA^1.$$ 
Let $U_1:=\AA^1\setminus \uo,$ $U_2:=\AA^1\setminus \vo$
and $S:=\AA^1\setminus \uo \ast \vo.$ Set
$$\tilde{f}(x,y):= \prod_{i=1}^{p} (x-x_i) \prod_{j=1}^q(y-x-y_j)
\prod_{i,j} (y-(x_i+y_j)) $$ and let $f\in \CC[x,y]$ be the
associated reduced polynomial. If the cardinality of 
$\uo \ast \vo$ is equal to $ p\cdot q,$ 
 we call $\uo \ast \vo$ {\em generic}.
Let
$ {\bf w}:=\{ (x,y)\in \AA^2 \mid f(x,y) =0\}$
and let $U:=\AA^2 \setminus {\bf w}.$ 
The space $U$ is equipped with the two projections
$$\pr_1: U \To U_1,\,  (x,y) \Mapsto x,\quad  \pr_2: U \To S,\, (x,y) 
\Mapsto y ,$$
and the {\it subtraction map}
$$ \q: U \To U_2,\quad  (x,y) \Mapsto y-x\,.$$
Let $ j:U \To \PP^1_{S}$ denote the canonical inclusion.
We further fix a basepoint $(x_0,y_0)$ of $U$ and let 
$U_0:=\pr_2^{-1}(y_0).$

\begin{defn}\label{defmidconv}{\rm  For  $\V_1 \in \LS_R(U_1)$ and
$\V_2 \in \LS_R(U_2)$ consider the tensor product
$$\V_1 \tim \V_2= \pr_1^*\V_1\otimes \q^*\V_2\quad \in\quad  \LS_R(U).$$
Then the {\em middle convolution} of $\V_1$ and
$\V_2$ is the higher direct image 
$$ \V_1 \ast \V_2 := R^1 (\opr_2)_*\left(j_*(\V_1\tim \V_2)\right).$$}
\end{defn}

The following proposition follows from the local triviality 
of $\pr_2$ (using the arguments of \cite{Hartshorne}, Prop.~8.1):

\begin{prop}\label{dimensione1}  The  middle convolution 
$\V_1\ast \V_2$ is a local system on $S$ whose stalk 
$(\V_1\ast \V_2)_{y_0}$ at $y_0$ is canonically isomorphic to 
the parabolic cohomology 
$H^1_p(U_0,\V_1\tim\V_2|_{U_0}).$ 
\end{prop}

\subsection{The basic topological setup}\label{basictop}
 Throughout this and the following subsections, we
assume   that the coefficienty ring $R$ is a field $K$ 
and that $\uo \ast \vo$
is  generic.
The first projection
yields an identification of $U_0=\pr_2^{-1}(y_0)$ with
$$\AA^1\setminus {\bf d},\quad {\rm where}\quad
{\bf d}=\{x_1,\ldots,x_p,y_0-y_1,\ldots,y_0 - y_q\}.$$
Using a suitable homeomorphism 
$\PP^1\to \PP^1$ we can assume
 that we are in the following situation: 
The sets
$\uo=\{x_1,\ldots,x_p\},\,\vo=\{y_1,\ldots,y_q\},\, \{y_0\}$ are
element-wise real and 
$$0<x_1<x_2<\ldots < x_p <
y_0-y_1<y_0-y_2<\ldots< y_0-y_q\,.$$ Moreover, we can assume that
\begin{equation}\label{differenceofvu}
|x_p-x_1|< |y_{i+1}-y_i|\quad {\rm for}\quad  i=1,\ldots,q-1. \end{equation}
 We choose
generators $\al_1,\ldots,\al_{p+q}$ of $\pi_1(\AA^1\setminus {\bf d},x_0)$ as follows:

%\pagebreak 
\vspace{.6cm}
\begin{center}
\includegraphics{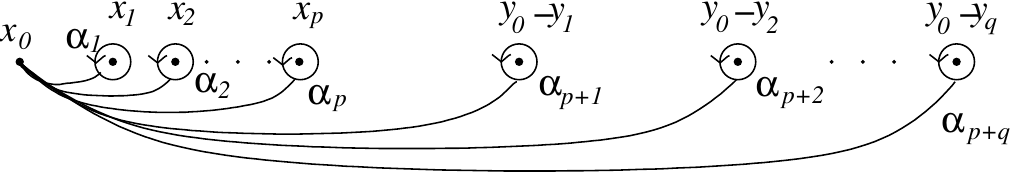}
\end{center}
%\caption{\label{zopfbild} The generators $\al_1,\ldots,\al_{p+q}$}

\noindent We also choose generators $\delta_{i,j},\,
i=1,\ldots,p,\,j=1,\ldots ,q$ of $\pi_1(S,y_0)$ as follows:
\vspace{.6cm}
\begin{center}
\includegraphics{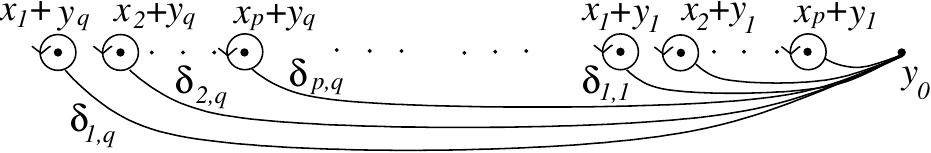}
\vspace{.6cm}
\end{center}
The long exact homotopy
sequence associated to the locally trivial 
fibration $\pr_2:U\to S$ yields a split short exact sequence
$$ 1\To\pi_1(\AA^1\setminus {\bf d},x_0)\To \pi_1(U,(x_0,y_0))\To  \pi_1(S,y_0)\To 1.$$
By embedding 
$S\setminus \{y_1+x_0,\ldots,y_q+x_0\}$
 and $\AA^1\setminus \d=U_0$ into $U$  we can think of the paths 
$\alpha_k$ and $\delta_{i,j}$ as loops in $U.$ 
 Then the path
${}^{\delta_{i,j}}\alpha_k:= \delta_{i,j}\alpha_k\delta_{i,j}^{-1}$
becomes a loop in $U$ and it is in the
kernel of the third arrow of the above exact sequence. Therefore it has a unique
inverse image under the second arrow; this defines ${}^{\delta_{i,j}}\alpha_k$ as
an element of the fundamental group of the fiber. This gives us an action of
$\pi_1(S,y_0)$ on $\pi_1(U_0,(x_0,y_0))$. It is standard, that for
$i=1,\ldots,p,$ the following formula holds (cf.~\cite{DR07}, Section~4.1):
\begin{equation}\label{zopf}
 ({}^{\delta_{i,1}}\alpha_1,\ldots,{}^{\delta_{i,1}}\alpha_{p+1})=
 (\alpha_1,\ldots, \alpha_{i-1},\alpha_i^{\alpha_{p+1}},
\alpha_{i+1}^{[\alpha_i,\alpha_{p+1}]},
 \ldots,\alpha_p^{[\alpha_i,\alpha_{p+1}]}, \alpha_{p+1}^{\alpha_i
\alpha_{p+1}}), \end{equation}
where $[\alpha_i,\alpha_j]=\alpha_i^{-1}\alpha_j^{-1}\alpha_i\alpha_j$ 
and $$\alpha_i^{\alpha_{p+1}}=\alpha_{p+1}^{-1}\alpha_i\alpha_{p+1},\quad \alpha_{k+1}^{[\alpha_i,\alpha_{p1}]}=([\alpha_i,\alpha_{p+1}])^{-1}
\alpha_{k+1} [\alpha_i,\alpha_{p+1}]\quad {\rm etc.}$$

\subsection{The rank of the middle convolution}
Let $\V_1 \in \LS_R(U_1)$ and $\V_2 \in \LS_R(U_2).$ The homotopy base 
$\al_1,\ldots,\al_{p+q}$ on $U_0$ given in Section~\ref{basictop}
 induces homotopy bases 
on $U_1$ and $U_2$ by looking at the images of the maps
$\pr_1|_{U_0}$ and $\d|_{U_0}$ (resp.).  Let 
the   monodromy tuples of $\V_1$ and $\V_2$ with respect
to these homotopy bases be 
$$T_{\V_1}=(A_1,\ldots,A_{p+1})\in \GL(V_1)^{p+1}$$ and
$$T_{\V_2}=(B_1,\ldots,B_{q+1})\in \GL(V_2)^{q+1}$$
(resp.).  Then the monodromy
tuple of $\V_1\tim \V_2|_{U_0}$ is of the following form:
\begin{multline} \label{Tcirc}
T_{\V_1\tim \V_2|_{U_0}} = (C_1= A_1\otimes 1_{V_2}\,,\,\ldots,\,
C_p= A_p\otimes 1_{V_2}\,,\\
 C_{p+1}=1_{V_1}\otimes B_1\,,\,\ldots\,,\, C_{p+q}=1_{V_1}\otimes 
B_q\,,\,C_{p+q+1}=A_{p+1}\otimes B_{q+1})\in \GL(V_1\otimes V_2)^{p+q+1}\,.
\end{multline}

\begin{prop}\label{dimensione}  Suppose that $R=K$ is a field and
that $\V_1\in \LS_K(U_1)$ has no 
global sections.  
Let $\dim_K V_i=n_i,\,i=1,2.$ Then
\begin{multline}\label{dimensionconv} \rk(\V_1\ast \V_2)= (p+q-1)n_1 n_2 -
\sum_{i=1}^p n_2\dim_K
\ker(A_i - 1_{V_1}) \\
-\sum_{j=1}^q n_1\dim_K \ker(B_j - 1_{V_2})
-\dim_K\ker(A_{p+1}\otimes B_{q+1}-1_{V_1\otimes
V_2})\,.\end{multline}
\end{prop}

\proof We have to determine the dimension
of the stalk 
$(\V_1\ast \V_2)_{y_0}$ which is isomorphic to  the 
parabolic cohomology group
$H^1_p(U_0,\V_1\tim\V_2|_{U_0})$  by  Prop.~\ref{dimensione1}.
It follows from \eqref{Tcirc} and the properties of the
tensor product that
$$\dim_K \ker(C_i-1_{V_1\otimes V_2})= n_2\dim_K
\ker(A_i - 1_{V_1}),\quad {\rm for} \quad i=1,\ldots, p\,,$$ and
$$\dim_K \ker(C_i-1_{V_1\otimes V_2})= n_1\dim_K
\ker(B_{i-p} - 1_{V_2}),\quad {\rm for} \quad i=p+1,\ldots, p+q\,.$$ 
Let $\V=\V_1\tim \V_2|_{U_0}.$ 
 The Euler-Poincar\'e Formula implies that 
\begin{eqnarray}
\chi(\PP^1,j_*\V)&=&h^0(\PP^1,j_*\V)
-h^1(\PP^1,j_*\V)+h^2(\PP^1,j_*\V)
\nonumber \\
&=& \chi(U_0)\cdot \rk(\V) +\sum_{i=1}^{p+q} 
\dim(\ker(C_i-1)) + \dim(\ker(C_{p+q+1}-1))\nonumber \\
&=& (1-p-q) n_1n_2+ \sum_{i=1}^p n_2\dim_K
\ker(A_i - 1_{V_1}) \nonumber \\
&&\quad 
+\sum_{j=1}^q n_1\dim_K \ker(B_j - 1_{V_2})
+\dim_K\ker(A_{p+1}\otimes B_{q+1}-1_{V_1\otimes
V_2}).\nonumber \end{eqnarray}
Since  $\V_1$ has no global sections, the local system
$\V$ has no global sections. Therefore, 
$$0=h^0(U_0,\V)=h^0(\PP^1,j_*\V)=h^2(\PP^1,j_*\V),$$
which implies 
the claim (the third equality follows from Poincar\'e duality, using the
same arguments for the dual of $\V$). 
 \Endproof

\subsection{The twisted 
evaluation map}\label{basictop1}

As in the last section, let $$\rho:=\rho_{\V_1\tim \V_2|_{U_0}}:\pi_1(S,y_0)\To \GL(V_1\otimes V_2)$$
be the  monodromy representation of $\V_1\tim \V_2|_{U_0}.$
The Hochschild-Serre spectral sequence implies  that 
$H^1(U_0,\V_1\tim \V_2|_{U_0})\simeq H^1(\pi_1(U_0),V_1\otimes V_2)\,.$
The group cohomology $H^1(\pi_1(U_0),V_1\otimes V_2)$ is the 
quotient of
$C^1(\pi_1(U_0),V_1\otimes V_2)$ by $B^1(\pi_1(U_0),V_1\otimes V_2),$ where 
\begin{eqnarray}\label{eqcoc}
C^1(\pi_1(U_0),V_1\otimes V_2)&:=&
\{(\delta:\pi_1(U_0)\to V_1\otimes V_2)\quad \mid\nonumber \\
&&\quad \quad  \delta(\alpha \beta)=
\delta(\alpha) \rho(\beta)+\delta(\beta)
\quad \forall\; \alpha,\beta\in \pi_1(U_0)\}\end{eqnarray}
and 
$$B^1(\pi_1(U_0),V_1\otimes V_2):=\{ \delta_v\mid v \in V_1\otimes V_2,\quad 
\delta_v(\gamma)=v(1-\rho(\gamma))\quad \forall\; \gamma\in \pi_1(U_0)\}.$$

\begin{defn}{\rm The linear map
$$
 \tau :C^1(\pi_1(U_0),V_1\otimes V_2)\To (V_1\otimes V_2)^{pq}\, ,$$
$$\delta \Mapsto 
\left(\,\delta([\alpha_1,\alpha_{p+1}]),\cdots,\,
\delta([\alpha_p,\alpha_{p+1}]),\ldots,\,
\delta([\alpha_1,\alpha_{p+q}]),\cdots,
\,\delta([\alpha_p,\alpha_{p+q}])\,\right),$$
is called the {\em twisted evaluation map}.}
\end{defn}

\begin{lem}\label{3.4}  Assume that $\V_1\in \T_K(U_1)$ 
(see Def.~\ref{defconvosheaf}) with 
$T_{\V_1}=(A_1,\ldots, A_{p+1}),$
and assume that $\V_2\in \LS_K(U_2)$ is a rank one system 
 with $$T_{\V_2}=(\lambda_1,\ldots ,\lambda_{q+1})\in \GL_1(K)^{q+1}
\quad {and} \quad \lambda_i\not=1\quad {for} \quad i=1,\ldots,q.$$
Then the  kernel of the twisted evaluation map 
coincides with the coboundaries 
 $B^1(\pi_1(U_0),V),$ where 
$V=V_1\otimes V_2\simeq V_1.$ 
 Hence, the 
twisted evaluation map $\tau$   induces an embedding
$$ \tau :H^1(\pi_1(U_0),V_1\otimes V_2) \To (V)^{pq}.$$ 
\end{lem}

 \proof The cocycle relation \eqref{eqcoc} implies
\begin{equation}\label{eiei1}\delta([\alpha_i,\alpha_{p+j}])= \delta(\alpha_i)(\lambda_j-1)+
\delta(\alpha_{p+j})(1-A_i),\quad {\rm for}\quad 
i=1,\ldots,p,\,j=1,\ldots,q \,.\end{equation} 
By definition, the coboundaries $B^1(\pi_1(U_0),V)$ are therefore contained 
in  $\ker(\tau).$
Conversely, assume that $\delta\in \ker(\tau).$
Since $\lambda_j$ is assumed to be $\not=1$ one has
\begin{equation}\label{eqcotz} \delta(\alpha_i)=\frac{1}{1-\lambda_j}\delta(\alpha_{p+j})(1-A_i),
\quad {\rm for}\quad  i=1,\ldots,p,\,j=1,\ldots,q\,.\end{equation}
Set $v_j:=\frac{1}{1-\lambda_j}\delta(\alpha_{p+j}).$ By \eqref{eqcotz},
$$ v_j(1-A_i)=v_{j'}(1-A_i)\quad {\rm for} \quad 
i=1,\ldots,p,\quad {\rm and}\quad j,j'=1,\ldots,q \,.$$
If $v_j\not=v_{j'}$ for some $j\not= j',$ then the above equality shows that 
the vector $v_j-v_{j'}$ spans a trivial 
$\langle A_1,\cdots, A_p\rangle$-submodule of $V_1.$ But 
since $\V_1$ was assumed
to be contained in $\T_K(U_1)$ this is impossible. Thus 
$v_j=v_{j'}$ for $j\not= j'$ and 
$$\delta(\alpha_i)=\frac{1}{1-\lambda_q}\delta(\alpha_{p+q})(1-\rho(\alpha_i))\,.$$ Consequently, the element $\delta$ is an element of $B^1(\pi_1(U_0),V).$
\Endproof 

\subsection{Convolution with Kummer sheaves.}\label{secirrkumm1}
We work in the setup of the last sections with $U_1=\AA^1\setminus 
\uo\,\, (\uo=\{x_1,\ldots,x_p\})$ and $U_2=\GG_m=\AA^1\setminus \{0\}.$
 Let $K$ be a field and let $\V_\chi\in \LS_K(\GG_m)$ be the {\it Kummer sheaf} 
associated to a character
$$\chi:\pi_1(\GG_m)\To K^\times,\quad \alpha \Mapsto \lambda\not=1,$$
where $\alpha $ is a generator of $ \pi_1(\GG_m)$ moving counterclockwise
around the origin.  Let
 $\V\in \T_K(U_1).$ 
By Prop.~\ref{dimensione1}, the middle convolution $\V\ast \V_\chi$ is a local system 
on $S=\AA^1_y\setminus \uo,$ and can therefore be seen as 
a local system on $U_1.$  In this way, we obtain a functor 
$$ \MC_\chi: \, \LS_K(U_1) \To \LS_K(U_1),\quad \V\Mapsto \MC_\chi(\V):=\V \ast \V_\chi.$$

The middle convolution $\MC_\chi(\V)$ is a sub-local system of 
$R^1(\pr_{2})_*(\V \tim\V_\chi).$ 
We have 
$$R^1(\pr_{2})_*(\V \tim \V_\chi)_{y_0}\simeq H^1(U_0,\V\tim\V_\chi |_{U_0})
\simeq H^1(\pi_1(U_0,(x_0,y_0)),V),$$ where $V=V_1\otimes V_2\simeq V_1$ is the stalk of  $\V\tim\V_\chi |_{U_0}$
at $(x_0,y_0)$ with its natural action of $\pi_1(U_0,(x_0,y_0)).$
The  monodromy operation of 
$\delta\in \pi_1(S,y_0)$ on the stalk
$R^1(\pr_{2})_*(\V \tim \V_\chi)_{y_0}$ 
is given by conjugating the 
arguments of the  $1$-cocycles in $C^1(\pi_1(U_0),V)$
with the inverse of $\delta,$  cf.~\cite{dw03},~Lemma 3.3. 
Under the 
twisted evaluation map $\tau:C^1(\pi_1(U_0),V)
\to V^p,$ this action of $ \pi_1(S,y_0)$ on $H^1(\pi_1(U_0),V)$ 
defines an action of $\pi_1(S,y_0)$ on $V^p.$  Let
$\tilde{D}_{i}\, (i=1,\ldots,p)\in \GL(V^p)$ be the automorphisms
 induced in this way by the homotopy generators 
$\delta_i:=\delta_{i,1}$ of $\pi_1(S)$ ($i=1,\ldots,p$), where the $\delta_{i,1}$
are as 
in Section~\ref{basictop}.

\begin{lem}\label{lempochmat} Let $(A_1,\ldots,A_{p+1})\in \GL(V)^{p+1}$
be the monodromy tuple of $\V.$
Then for $i=1,\ldots,p,$ 
the automorphism $\tilde{D}_{i}\in GL(V^p)$ is a block matrix 
which is the identity block matrix (the block structure induced by 
$V^p$) 
outside the 
$i$-th block row and the $i$-th block row is as follows:
\[ 
              ( \lambda (A_1-1) , \ldots, \lambda(A_{i-1}-1)  , \lambda A_{i} , (A_{i+1}-1) ,\ldots, 
   (A_p-1)).\] 
\end{lem}

\proof The cocycle relation \eqref{eqcoc} implies that
\begin{equation}\label{eqcoomraus}
 \delta(\alpha[\beta,\gamma]\epsilon)=\delta(\alpha\epsilon)+
\delta([\beta,\gamma])\rho(\epsilon)\end{equation}
and $\delta(\alpha^{-1})=-\delta(\alpha)\rho(\alpha)^{-1}.$
If $m<i,$ then these formulas together 
with \eqref{zopf} yield
\begin{eqnarray}\delta[\alpha_m,\alpha_{p+1}]\tilde{D}_{i}&=&
\delta[\alpha_m^{(\delta_{i})^{-1}},
\alpha_{p+1}^{(\delta_{i})^{-1}}]\nonumber \\
&=& \delta[\alpha_m,\alpha_{p+1}^{\alpha_i\alpha_{p+1}}]\nonumber \\
&=&\delta(\alpha_m^{-1}\alpha_{p+1}^{-1}[\alpha_i,\alpha_{p+1}]\alpha_m
[\alpha_{p+1},\alpha_i] \alpha_{p+1})\nonumber \\
&=&\delta[\alpha_i,\alpha_{p+1}]\lambda_1(A_m-1)+\delta[\alpha_m,\alpha_{p+1}]\,.
\nonumber \end{eqnarray}
If $m=i$ then 
\begin{eqnarray}\delta[\alpha_i,\alpha_{p+1}]\tilde{D}_{i}&=&
 \delta[\alpha_i^{\alpha_{p+1}},
\alpha_{p+1}^{\alpha_i\alpha_{p+1}}]\nonumber \\
&=&\delta(\alpha_{p+1}^{-1}\alpha_i^{-1}[\alpha_i,\alpha_{p+1}]\alpha_i\alpha_{p+1})\nonumber \\
&=&\delta[\alpha_i,\alpha_{p+1}]A_i\lambda_1\nonumber\,, 
\end{eqnarray}
where in the last equation we have used \eqref{eqcoomraus} and  
$\delta(1)=0.$
If $m>i,$ then 
\begin{eqnarray}\delta[\alpha_m,\alpha_{p+1}]\tilde{D}_{i}&=&
 \delta[\alpha_m^{[\alpha_i,\alpha_{p+1}]},\alpha_{p+1}^{\alpha_i\alpha_{p+1}}]\nonumber \\
&=&\delta[\alpha_i,\alpha_{p+1}](A_m-1)+\delta[\alpha_m,\alpha_{p+1}]\,
\nonumber. \end{eqnarray}
This proves the claim.\Endproof

Let 
$$ {\cal K}=\{\;(w_1,\ldots,w_p)\mid w_i\in \im(A_i-1)\;\}\subseteq V^p,$$ 
$$ {\cal L}=\{\;(w_1A_2\cdots A_p,w_2A_3\cdots A_p,\ldots,w_p)\mid w_i\in \im(A_1\cdots A_p\lambda_1-1)\;\}\subseteq V^p\,,$$
and let $W:={\cal K} \cap {\cal L}\subseteq V^p.$
The automorphisms $\tilde{D}_i\, (i=1,\ldots,p)$ can easily be seen to stabilize the subspaces $\K$ and $\L$ of $V^p$ and hence
$W.$ Let ${\bf A}=(A_1,\ldots,A_{p+1}),$ and let 
$$\MC_\lambda({\bf A}):=(D_1,\ldots,D_{p+1})\in \GL(W)^{p+1}$$ 
be the tuple of linear 
transformations on $W$ induced by $\tilde{D}_1,
\ldots,\tilde{D}_p$ and $\tilde{D}_{p+1}:=(D_1\cdots D_p)^{-1}$ (resp.).

\begin{thm}\label{remlab} Let $\V\in \T_K(U_1)$ and 
let ${\bf A}\in \GL(V)^{p+1}$ be the monodromy 
tuple of $\V$ with respect to 
$\alpha_1,\ldots,\alpha_p.$ 
Then the stalk of $\MC_\chi(\V)$ at $y_0$ is canonically isomorphic to $W$
 and the monodromy tuple of $\MC_\chi(\V)$ with respect to 
the homotopy base $\delta_1,\ldots,\delta_{p}$ 
is given by $\MC_\lambda({\bf A})\in \GL(W)^{p+1}.$ 
\end{thm}

\proof The parabolic cohomology 
$H^1_p(U_0,\V\tim \V_\chi|_{U_0})=\MC_\chi(\V)_{y_0}$ can be seen 
as the subspace of 
$H^1(\pi_1(U_0),V)$ consisting of the elements $[\delta]\in 
H^1(\pi_1(U_0),V),$ where $\delta\in C^1(\pi_1(U_0),V)$ is a
{\it parabolic cocycle}, i.e., 
\begin{equation}\label{eiei}\delta(\gamma)\in \im(\rho_{\V\tim \V_\chi|_{U_0}}(\gamma)-1),\quad \forall \gamma \in \pi_1(U_0)\,,\end{equation} 
see \cite{dw03}, Lemma~1.2. Let $W'\subseteq V^p$ 
be the image of the parabolic cohomology group 
$H^1_p(U_0,\V\tim \V_\chi|_{U_0})$ in $V^{p}$ 
under  the twisted evaluation map and let $(w_1,\ldots,w_p)$ $\in W'.$ 
We have $w_i=\delta([\alpha_i,\alpha_{p+1}]),$ where $\delta$ is some parabolic 
cocycle and 
$\alpha_i,\,i=1,\ldots,p+1,$ is the homotopy base of 
$\pi_1(U_0)$ chosen according to the basic topological setup. 
By \eqref{eiei} we have 
 $\delta(\alpha_i)\in \im(A_i-1).$ The formula in \eqref{eiei1} implies therefore that 
$w_i\in \im(A_i-1)$ for $i=1,\ldots,p,$ so  $W'\subseteq \K.$ Similarly, the condition
 $$\delta(\alpha_1\cdots\alpha_{p+1})\in 
\im(A_1\cdots A_p \cdot \lambda_1-1) \textrm{ for all }  [\delta] \in H^1_p(U_0,V)\,$$ forces $W'\subseteq \L$ and 
hence  $W'\subseteq \K\cap \L.$ By Prop.~\ref{dimensione},
the rank of $\MC_\chi(\V)$ is equal to 
$$ pn- \sum_{i=1}^p \dim(\ker(A_i-1))-\dim(\ker(A_{p+1}\lambda^{-1}-1)).$$
Let $V^*$ denote the dual space of $V$ and consider the column space
$(V^*)^p$ as dual space of $V^p.$ Let 
$\K^*\leq (V^*)^p$ and $\L^*\leq (V^*)^p$ be the annihilators of the spaces 
$\K$ and $\L.$ It is not hard to verify that since $\lambda\not=1,$ 
one has $\K^*\cap \L^*=0,$ see \cite{dr00}, Lemma~2.7. 
This implies that 
$$\dim(W)=\dim(\K\cap \L)=\dim({\rm Ann}_{V^p}(\K^*\oplus \L^*))=
pn-\dim(\K^*)-\dim(\L^*),$$ compare to 
  \cite{dr00}, Lemma~2.7. It follows therefore from dimension reasons 
that 
$W'=W,$ which implies the claim.  
 \Endproof

\begin{rem} One can switch from the homotopy base $\delta_1,\ldots,\delta_p$
of $\pi_1(S,y_0)$ back to the original homotopy base $\alpha_1,\ldots,
\alpha_p$ of $\pi_1(S,x_0)$
by connecting $x_0$ with $y_0$ using a path $\gamma$ in the lower half-plane. 
In this way, one obtains a well defined transformation $\MC_\chi$ on the level
of representations of $\pi_1(S,x_0).$ 
\end{rem}

By construction, the tuple transformation $\MC_\lambda$ considered 
above is the dual of the tuple transformation $\MC_\lambda$ 
considered in \cite{dr03} and 
\cite{dr00} (we use right-actions on row spaces, where
\cite{dr03} and 
\cite{dr00} use left-actions on column spaces). Since irreduciblity is 
preserved under dualization, Theorem~2.4 of \cite{dr03} implies the following:

\begin{cor}\label{cor2.4} If $\V\in \T_K(U_1),$ then the
local system $\MC_\chi(\V)$ is again 
irreducible.
\end{cor}

\subsection{Irreducibility criteria for  the middle convolution.}\label{secirrkumm}

\begin{thm}\label{thmirrd} Let $\V_1\in \T_K(U_1)$ with
$T_{\V_1}=(A_1,\ldots, A_{p+1})\in \GL(V)^{p+1}$
and assume that $\V_2\in \LS_K(U_2)$ is a rank-one system 
 with $$T_{\V_2}=(\lambda_1,\ldots ,\lambda_{q+1})\in \GL_1(K)^{q+1}
\quad {\rm and} \quad \lambda_i\not=1\quad {\rm for} \quad i=1,\ldots,q.$$
Assume that $\uo\ast \vo$ is generic.
Then the local system $\V_1\ast\V_2\in \LS_K(S)$ is irreducible if
$$ (p-2)n-
\sum_{i=1}^p\dim_K(\ker(A_i-1))>0\,.$$
\end{thm}

\proof This follows from  induction on $q.$ For $q=1$ this is 
Cor.~\ref{cor2.4}. If $q>1,$ then we can assume that 
${\V_1}\ast \tilde{\V}_2$ is irreducible, where 
$$ \tilde{\V}_2\in \LS_K(\AA^1\setminus \{y_2,\ldots,y_q\})\quad {\rm with}\quad T_{\tilde{\V}_2}=(\lambda_2,\ldots,\lambda_q,\lambda_1\cdot \lambda_{q+1})\,.$$
Let 
$$ p_1:(V^p)^q\To V^p ,\, (v_1,\ldots,v_q)\Mapsto v_1\,,$$
and 
$$p_2:(V^p)^q\To (V^p)^{q-1} ,\, (v_1,\ldots,v_q)\Mapsto (v_2,\ldots,v_q)\,.$$
Let $$G_1:=\langle {\delta_{i,1}},\, i=1,\ldots,p\rangle\,\leq \,\pi_1(S,y_0)$$ and 
$$G_2:=\langle {\delta_{i,j}},\, i=1,\ldots,p,\, j=2,\ldots,q\rangle\,\leq \,\pi_1(S,y_0)\,.$$
By Lemma~\ref{3.4}, there is an embedding
$$\tau:(\V_1\ast \V_2)_{y_0}=H^1_p(U_0,\V_1\tim \V_2|_{U_0})\To (V^p)^q.$$
We can further assume that $y_1=0.$ 
 By the assumptions defined in the basic topological setup,
the $G_1$-module
$\im(p_1\circ \tau)$ is then isomorphic to the 
monodromy representation of 
$\MC_{\lambda_1}(\V_1)$ and is hence irreducible
by Cor.~\ref{cor2.4}.  Since  $\pi_1(S,y_0)$ is the free product 
of $G_1$ and $G_2,$ and since 
the rank of $\MC_{\lambda_1}(\V_1)\geq pn-\sum_{i=1}^p \dim(\ker(A_i-1))$ 
by Prop.~\ref{dimensione}, 
 it follows
that 
$(\V_1\ast \V_2)_{y_0}$ contains an irreducible 
$\pi_1(S,y_0)$-submodule 
$W_1$ of
rank greater than or equal to 
$$n_1=pn-\sum_{i=1}^p \dim(\ker(A_i-1))\,.$$ 

By the induction hypothesis, the $G_2$-module
$$ (\V_1\ast \tilde{\V}_2)_{y_0}=
H^1_p(\AA^1\setminus (\uo\cup \{y_0-y_2,\ldots,y_0-y_q\}),
\V_1\tim\tilde{\V}_2|_{\AA^1\setminus (\uo\cup \{y_0-y_2,\ldots,y_0-y_q\}})$$ is irreducible. Moreover, the image of 
 the $G_2$-module-homomorphism 
 $$p_2\circ \tau: (\V_1\ast \V_2)_{y_0} \To 
V^{p(q-1)}$$ coincdes with the image of of 
 the $G_2$-module-homomorphism 
$$ \tau: (\V_1\ast \tilde{\V}_2)_{y_0}
\To V^{p(q-1)}.$$  By Prop.~\ref{dimensione}, the $\pi_1(S,y_0)$-module
$(\V_1\ast \V_2)_{y_0}$ contains therefore an irreducible 
submodule 
 $W_2$ of rank greater than or equal to $$ n_2:=(p+q-3)n-
\sum_{i=1}^p\dim_K(\ker(A_i-1))\,.$$ 
By Prop.~\ref{dimensione}, the rank of $\V_1\ast \V_2$ is  smaller than or equal to
$$ n_3=(p+q-1)n-\sum_{i=1}^{p}\dim_K(\ker(A_i-1))\,.$$
One has 
$$n_1+n_2=2pn +qn-3n-2\cdot \sum_{i=1}^p\dim_K(\ker(A_i-1))$$
Thus 
$$ n_1+n_2-n_3\geq (p-2)n-
\sum_{i=1}^p\dim_K(\ker(A_i-1))\,.$$ 
By assumption,  $(p-2)n-\sum_{i=1}^p\dim_K(\ker(A_i-1))>0,$ thus
$ n_1+n_2>n_3.$ It follows that the intersection of the irreducible submodules 
$W_1$ and $W_2$ is 
non-trivial and hence $\V_1\ast \V_2$ is irreducible (the spaces
$W_1$ and $W_2$ clearly generate $(\V_1\ast \V_2)_{y_0}$ as a vector space).
 \Endproof

\subsection{The local monodromy of the middle convolution.}\label{secmonodrconv}

In the following, an expression $J(\alpha,l)\, (\alpha\in K,\, l\in \NN)$
denotes a Jordan block of length $l$ and eigenvalue $\alpha.$  

\begin{lem} \label{lemmonodromy1} Let $K$ be algebraically closed.
Let $\V_1\in \T_K(U_1)$ with
$T_{\V_1}=(A_1,\ldots, A_{p+1})\in \GL(V)^{p+1}$
and assume that $\V_2\in \LS_K(U_2)$ is a rank-one system 
 with $T_{\V_2}=(\lambda_1,\ldots ,\lambda_{q+1})\in \GL_1(K)^{q+1}$ and
$ \lambda_i\not=1\quad$ for $i=1,\ldots,q.$
Assume that $\uo\ast \vo$ is generic. Let $D_{i,j}:=\rho_{\V_1\ast \V_2}(\delta_{i,j}),$ where the $\delta_{i,j}$ are the generators of
$\pi_1(S)$ as above.  Then the following
holds:
Every $\lambda_j$  and every Jordan
block $J(\alpha,l)\not= J(1,1)$ occurring in the Jordan
decomposition of $A_i$ contribute a Jordan block $J(\alpha
\lambda_j,l')$ to the Jordan decomposition of $D_{i,j},$ where
$$ l':\;=\quad
  \begin{cases}
    \quad l &
                              \quad\text{\rm if $\alpha \not= 1,\lambda_j^{-1}$,} \\
    \quad  l-1& \quad \text{\rm if $\alpha =1$,} \\
    \quad l+1 & \quad \text{\rm if $\alpha =\lambda_j^{-1}$.}
  \end{cases}
  $$
  The only other Jordan blocks which occur in the Jordan
  decomposition of $D_{i,j}$ are blocks of the form $J(1,1).$
\end{lem}

\proof We can assume  that $y_j=0.$ Let
$S^\circ:=\{y\in S\mid 0<|y|<x_p+\epsilon\}\, ,$ where
$\epsilon\in \RR_+$ is very small.  The restriction of
$\pr_2$ to $\pr_2^{-1}(S^\circ)$ gives then rise to a variation
of parabolic cohomology groups which is essentially 
$\MC_{\chi}(\V_1),$ where $\chi$ is the character 
defined by sending a counterclockwise 
generator of $\pi_1(\GG_m)$ to $\lambda=\lambda_j$
(by the conventions defined in the basic topological
setup, the other components $\lambda_k\not= \lambda_k$ do not 
contribute to the variation, although they contribute to the 
rank of $\V_1\ast \V_2$).  
 Let $D_i$ be the $i$-th component of the monodromy tuple 
$\MC_\lambda(T_{\V_1})$ of $\MC_\chi(\V_1),$ compare to Thm.~\ref{remlab}.  
 Then $D_{i,j}$ is conjugated 
to a matrix of the form 
$$D_i \oplus \left(\bigoplus_{k=1}^{t}\J(1,1)\right),\textrm{ where }t=\rk(\V_1\ast \V_2)-\rk(\MC_\lambda(\V_1).$$ Therefore it suffices to 
prove the claim for $D_i$ instead of $D_{i,j}.$ 
 It is not hard
to see that one has an isomorphism 
\begin{equation}\label{eq4.12}\phi:\im({D}_{i}-1) \To \im (A_i-1),
\quad (0,\ldots,w,0,\ldots,0)\Mapsto w(A_i-1),
\end{equation}
such that $\phi\circ {D}_{i} = \lambda A_i\circ \phi$
(use \cite{dr00}, Rem.~3.1, and the same arguments as \cite{dr00}, Lemma~4.1).
The claim follows from this by an induction on the number 
of Jordan blocks.
\Endproof

The following lemma gives the monodromy at infinity in the  case
of the convolution with Kummer sheaves:

\begin{lem}\label{lemmonodromy2} Let 
$\V\in \T_K(U_1)$  with 
$T_{\V}=(A_1,\ldots,A_{p+1})\in \GL(V)^{p+1}$ and let 
$\V_\chi$ be a Kummer sheaf associated to a non-trivial 
character 
$\chi: \pi_1(\GG_m(\CC))\to K^\times,\, \gamma \mapsto 
\lambda\,.$ 
 Let
$(D_1,\ldots,D_{p+1})$ be the monodromy tuple 
of $\MC_\chi(\V).$  Then the
following holds:
Every Jordan
block $J(\alpha,l)$ occurring in the Jordan
decomposition of $A_{p+1}$ contributes a Jordan block $J(\alpha
\lambda^{-1},l')$ to the Jordan decomposition of $D_{p+1},$ where
$$ l':\;=\quad
  \begin{cases}
    \quad l, &
                              \quad\text{\rm if $\alpha \not= 1,\lambda$,} \\
    \quad  l-1& \quad \text{\rm if $\alpha =\lambda$,} \\
    \quad l+1 & \quad \text{\rm if $\alpha =1$.}
  \end{cases}
  $$
  The only other Jordan blocks which occur in the Jordan
  decomposition of $D_{p+1}$ are blocks of the form $J(\lambda^{-1},1).$
\end{lem}

\proof The claim follows  
from Thm.~\ref{remlab} and \cite{dr00}, Lemma~4.1 (b), using an induction on the 
number of Jordan blocks.\Endproof

\section{Convolution of  \'etale sheaves} \label{secetaleconN}

\subsection{Etale sheaves.}\label{secetaledef}

Let $R$ be either a field or a normal integral domain which is 
of finite type over $\ZZ$ and let $S$ be an
regular integral  scheme of finite type 
over $R.$ If $s$ is a closed point of $S,$ then $\bar{s}$ denotes a geometric point
which extends $s.$ 
If $U$ is a scheme over $S$ 
and if  $\phi:S'\to S$ is  a morphism of schemes, then the basechange
of $U$ which is defined by $\phi$  is denoted by $U_{S'}.$

Let $\ell$ be an invertible prime in 
$R,$ and let 
$F$ be either a subfield of $\bar{\FF}_\ell,$
or a finite extension of $\QQ_\ell,$ or $\bar{\QQ}_\ell.$ 
The category of constructible $F$-sheaves
on $S$ is denoted by $\Constr_F(S).$ 
If $\pi:X\to S$ is smooth and if $\V\in \Constr(X),$ then
one has the notion of the higher direct  image $R^i\pi_*(\V)\in 
\Constr(S)$ and the notion of the higher direct image with
compact support $R^i\pi_!(\V)\in 
\Constr(S),$ see \cite{MilneEC}. 

The category of constructible $F$-sheaves on $S$ which are
lisse in the sense 
of \cite{SGA5}
is denoted by $\LS_F^\et(S).$  Any $\V\in \LS_F^\et(S)$
corresponds to its {\it monodromy representation}
$$\rho_\V:\pi_1(S,\bar{\eta})\to \GL(\V_{\bar{\eta}}),$$
see  \cite{SGA5}, V, VI. 
An object $\V\in \LS_F^\et(S)$ is called 
{\it tamely ramified} if $\rho_\V$ factors over the 
tame fundamental group $\pi_1^\tame(S,\bar{\eta}).$ 
The category of 
tamely ramified local systems is denoted by $\LS_F^\tame(S).$  
We call $\V\in\LS_F^\tame(S)$ {\it geometrically (absolutely) irreducible},
if $\rho_\V|_{\pi_1(S,\bar{s})}$ is (absolutely) irreducible 
for some (and hence for every, see \cite{Katz90}, 8.17.13)
 geometric point $\bar{s}$ of $S.$ 

 Let $F$ be either a finite extension of $\QQ_\ell$ or $\bar{\QQ}_\ell,$
and let $S/R$ as above, where $R$ is finitely generated over 
$\ZZ.$ 
A sheaf  $\V\in \Constr_F(S)$ is called
{\it pure of weight $k\in \ZZ$} if for any 
closed point $s$ of $S,$  the eigenvalues $\alpha\in \bar{F}^\times$  
of the  geometric Frobenius element
$F_{s}$ under the natural  action of   $F_{s}$ on $\V_{\bar{s}}$
satisfy the following condition:
The image of $\alpha$ under any embedding $F\to \CC$ is an
 algebraic number which is  of complex absolute value ${\rm N}(s)^{k/2},$ see \cite{DeligneWeil2}. 
 A sheaf $\V\in \Constr_F(S)$ is called
{\it mixed of weight $\leq k$} if it is an iterated extension
of sheaves which are pure of weight $\leq k.$ If $\V\in \Constr_F(S)$
is mixed of weight $\leq k,$ then the largest quotient 
of $\V$ which is pure of weight $k$ is denoted by $W^k(\V).$

\subsection{The middle convolution of  \'etale sheaves .}\label{secvondefetale} 

In this section, let $R$ be a field or a normal  integral domain which is 
of finite type over $\ZZ.$
 Let $\uo,\vo\subset \AA^1_R=\Spec(R[x])$ be reduced subschemes, defined by 
the vanishing of the polynomials
$$g_i(x)=\prod_j(x-x_{j_i})\in R[x],\quad i=1,2 \quad \textrm{ (respectively),}$$ 
 with $x_{j_i}$ \'etale over $R.$
Let $g_1\ast g_2(x)\in R[x]$ be the reduced polynomial 
associated to  the polynomial 
$$\prod_{j_1,j_2}(x-(x_{j_1}+x_{j_2}))\in R[x]$$
 and let $ \uo \ast \vo \subset \AA^1_R$ be the subscheme defined 
by the vanishing of $g_1\ast g_2.$
Let $U_1:=\AA^1_R\setminus \uo,$
$U_2:=\AA^1_R\setminus \vo$ and $S:=\AA^1_R\setminus \uo \ast \vo.$
 Let further $$U:=
{\rm Spec}(R[x,y,\frac{1}{g_1(x)\cdot g_2(y-x)\cdot g_1\ast g_2(y)}]).$$ 
Define
\begin{eqnarray*}
\pr_1\,:\, U \To U_1, &\textrm{ by }& (x,y) \Mapsto x,\\
  \pr_2\,:\, U \To S,&\textrm{ by }& (x,y) \Mapsto y ,\\
 \q\,:\, U \To U_2,&\textrm{ by }& (x,y) \Mapsto y-x.\end{eqnarray*}
Let $j:U\to \PP^1_S=\PP^1\times S$ be the natural 
inclusion.
The second projection $\PP^1_S \to S$ is
denoted by $\opr_2.$

\begin{defn} \begin{enumerate}{\rm 
\item  The {\em middle convolution} of $\V_1 \in \LS_F^\et(U_1)$
and $\V_2 \in \LS_R^\et(U_2)$ is the constructible sheaf
$$ \V_1 \ast \V_2 := R^1 (\opr_2)_*(j_*(\V_1\tim \V_2))$$
on $S=\AA^1\setminus \uo\ast \vo,$
where $\V_1\tim \V_2:=\pr_1^*(\V_1)\otimes \q^*(\V_2)\,.$
\item The {\em naive-convolution} of $\V_1 \in \LS_F^\et(U_1)$ and
$\V_2 \in \LS_R^\et(U_2)$ is the constructible sheaf
$$ \V_1 \ast_\naive \V_2:=R^1 (\pr_2)_!(\V_1\tim \V_2) \in \Constr_F(S)\, .$$
\item 
 Let $\V \in \LS_{\bar{\QQ}_\ell}(U_1)$ and let $\L_\chi$ be the
lisse {\rm Kummer} sheaf on $U_2=\GG_m$ associated to a nontrivial
character $\chi:\pi_1^\tame(\GG_m)\to \bQl^\times.$ Then we 
set $$ \MC_\chi(\V):=\V\ast \L_\chi.$$}\end{enumerate}
\end{defn}

\subsection{First properties}
Assumptions and notation as in Section~\ref{secetaledef} and~\ref{secvondefetale}. 
\begin{prop}\label{deligweil}  Let $\V_i\in \LS_F^\tame(U_i),\,i=1,2.$
Then the following holds:
\begin{enumerate}
\item If $R$ is a subfield of $\CC,$ then 
the middle convolution $\V_1\ast \V_2$
is lisse on $\AA^1_R\setminus \uo_1\ast \uo_2.$
\item Let $R$ be a normal  integral domain which is 
of finite type over $\ZZ$ and which has a generic point
of characteristic zero.  Then the naive convolution 
$\V_1\ast_\naive \V_2$ and the middle convolution
$\V_1\ast \V_2$ are lisse and tame on $\AA^1_R\setminus \uo_1\ast \uo_2.$
\item If $R$ is an normal   integral domain which is 
of finite type over $\ZZ$ and if $\V_1$ and $\V_2$ 
 are pure of weight $n_1$ and $n_2$ (resp.), 
then 
$$\V_1\ast_\naive \V_2\quad \textrm{ is mixed of weights }\quad\leq n_1+n_2+1$$
and 
$$\V_1\ast \V_2=W^{n_1+n_2+1}\left(\V_1\ast_\naive V_2\right).$$
\end{enumerate}
\end{prop}

\proof   The first claim follows from 
\cite{dw03}, Thm. 3.2. Let $\F=\V_1\tim \V_2.$ Since the 
generic point of $R$ has characteristic zero, the sheaf 
$R^1 (\pr_2)_!(\F)$ is also tame by \cite{KatzSE}, 4.7.1 (i). 
It follows further from \cite{KatzSE}, 4.7.1 (ii), that 
$\V_1\ast_\naive \V_2=R^1 (\pr_2)_!(\F)$ is lisse on 
$\AA^1_R\setminus \uo_1\ast \uo_2.$  
The excision sequence yields an exact sequence 
of constructible sheaves:
\begin{equation}\label{exe1}
R^0(\opr_2|_D)_*((j_*\F)|_D)\to R^1(\pr_2)_!(\F) \to 
R^1(\opr_2)_*(j_*\F)\to 0.\end{equation}
 By \cite{KatzSE}, 4.7.1 (iii), the formation 
of $j_*\F$ and of $Rj_*\F$ commutes with arbitrary basechange on $S$  and 
$j_*\F|D$ is again lisse and also tame.
Hence $R^0(\opr_2|_D)_*((j_*\F)|_D)$ is lisse and tame 
on $\AA^1_R\setminus \uo_1\ast \uo_2.$
It follows therefore
 from \eqref{exe1} that the 
sheaf $R^1(\opr_2)_*(j_*\F)$ is also contained in 
$\LS_F^\tame(\AA^1_R\setminus \uo_1\ast \uo_2),$  proving~(ii).

 The sheaf  $\F=\V_1\tim \V_2$ on $U$ is mixed  of weight 
$\leq n_1+n_2,$ so  the sheaf $j_*\F$ is  
mixed of weights $\leq n_1+n_2$ by \cite{DeligneWeil2}, Corollaire~1.8.9.
This implies that also 
$R^0(\opr_2|_D)_*((j_*\F)|_D)$ is  
mixed of weights $\leq n_1+n_2.$
Again, by \cite{DeligneWeil2}, Thm.~3.3.1, the sheaf 
$R^1 (\pr_2)_!(\F)$ is mixed of weights $
\leq n_1+n_2+1.$  Moreover, by loc.~cit. Thm.~3.2.3,
the sheaf  $R^1(\opr_2)_*(j_*\F)$ is pure of weight $n_1+n_2+1.$ 
Thus, by \eqref{exe1}, the sheaf 
 $R^1(\opr_2)_*(j_*\F)$ is the weight-$n_1+n_2+1$-quotient
of $R^1(\pr_2)_!(\F).$
\Endproof

\begin{prop}\label{propverschwindibus} Let $\V_1 \in \LS_F^\tame(U_1)$ and 
$\V_2 \in \LS_F^\tame(U_2)$ 
and assume that either  $\V_1$ or $\V_2$ is  
 geometrically 
irreducible and nonconstant.  Then 
$$ R^i (\pr_2)_!(\V_1\tim \V_2) =0 \quad {\rm if}\quad  i\not=1$$
and 
$$ R^i (\opr_2)_*(j_*(\V_1\tim \V_2)) =0 \quad {\rm if}\quad  i\not=1.$$
 \end{prop}

\proof The map $\pr_2$ is affine and smooth of relative dimension
one. Since formation of $R^i (\pr_2)_!$ commutes with arbitrary 
base change by \cite{KatzSE},~4.7.1~(ii), proper basechange shows that  
$ R^i (\pr_2)_!(\V_1\tim \V_2) =0$ for $i\not=1,2.$
If $R^2 (\pr_2)_!(\V_1\tim \V_2) \not= 0,$ then
there exists a geometric point $\bar{s}$ of $s$ for which 
$$\left(R^2 (\pr_2)_!(\V_1\tim \V_2)\right)_{\bar{s}}
=H^2(U_{\bar{s}},(\V_1\tim \V_2)|_{U_{\bar{s}}})    \not=0.$$ 
Consider the maps
$ \pr_1=\pr_1|_{U_{\bar{s}}}: U_{\bar{s}}
 \to U_{1}$ and $\d=\d|_{U_{\bar{s}}}: U_{\bar{s}} \to U_{2}.$ 
By smooth basechange, 
$ (\V_1\tim \V_2)|_{U_{\bar{s}}} \simeq \pr_1^*(\V_1)\otimes \d^*(\V_2).$
The lisse sheaf $(\V_1\tim \V_2)|_{U_{\bar{s}}}$ 
corresponds thus 
to a tensor product representation 
$$ \rho_1\otimes \rho_2:\pi_1(U_{\bar{s}})\To \GL(V_1\otimes V_2),$$
where $\rho_i:\pi_1(U_{\bar{s}})\to \GL(V_i),\, i=1,2,$ 
is the monodromy 
representation associated to $\pr_1^*(\V_1),$ resp. $d^*(\V_2).$  
We have to show that the group of 
coinvariants of $\rho_1\otimes \rho_2$ vanishes. 
Assume that $\V_2$ is
 geometrically 
irreducible and nonconstant.
Let $G_2
\leq \pi_1(U_{\bar{s}})$ be the subgroup 
generated by the inertia subgroups at the points of $\vo.$ Since 
$\V_2$ is geometrically irreducible and since $\V_2$ is lisse at $\uo,$
the properties of the tensor product imply that 
the restriction of 
$\rho_1\otimes \rho_2$ to $G_2$ decomposes as a direct sum of copies 
of an irreducible and nontrivial representation of $G_2.$ 
Thus the group of coinvariants of $\rho_1\otimes \rho_2$
vanishes and $R^2 (\pr_2)_!(\V_1\tim \V_2) = 0.$ If $\V_1$ 
 geometrically 
irreducible and nonconstant, then repeat the above arguments for $\V_1.$ 

The map $\opr_2$ is  smooth projective of relative dimension
one. Proper basechange shows that  
$ R^i(\opr_2)_*(j_*(\V_1\tim \V_2)) =0$ for $i\not=0,1,2.$
The excision sequence yields an isomorphism 
$$R^2 (\pr_2)_!(\V_1\tim \V_2)\simeq R^2(\opr_2)_*(j_*(\V_1\tim \V_2)).$$ 
By what was proven above,
  $R^2(\opr_2)_*(j_*(\V_1\tim \V_2))=0.$ On the other hand, Poincar\'e
duality implies that $R^0(\opr_2)_*(j_*(\V_1\tim \V_2))$ is dual 
to $R^2(\opr_2)_*(j_*(\V_1^*\tim \V_2^*)),$ where $\V_i^*,\,i=1,2,$ denotes the
dual of $\V_i$ (resp.).
 Thus, by repeating the above arguments with $\V_i^*,\,i=1,2,$
one sees that $R^2(\opr_2)_*(j_*(\V_1^*\tim \V_2^*))=0$ and, consequently,
$R^0(\opr_2)_*(j_*\V_1\tim \V_2)=0.$
\Endproof

Suppose that $R$ 
is contained in $\CC.$  If  
$\V \in \LS_F^\et(U),\, U=\AA^1_R\setminus \uo,$ then 
the composition of the monodromy representation
$\rho_{\V}: \pi_1^\et(U,s)\to \GL_{n}(F)$ with 
the canonical homomorphism $ \pi_1(U(\CC),s)\to \pi_1^\et(U,s)$
defines a local system $\V^\an$ on $U(\CC).$ We call 
 $\V^\an$ the {\it analytification of} $\V.$ 

\begin{prop}\label{dimension} Suppose that $R\subseteq \CC.$ 
Let  $\V_i \in \LS_F^\et(U_i),\, U_i=\AA^1_R\setminus \uo_i,\,i=1,2$ and let
 $\V_i^\an \in \LS_F(U_i(\CC)),\,i=1,2,$ be the analytifications of 
$\V_i.$ 
Then 
$$  \V_1^\an*\V_2^\an=(\V_1*\V_2)^\an.$$ 
Especially, the monodromy representation 
 $\rho_{(\V_1*\V_2)^\an}$ of $(\V_1*\V_2)^\an$ is the composition of 
$\rho_{\V_1*\V_2}$ with the canonical homomorphism 
$\pi_1(S(\CC))\to \pi_1^\et(S).$ 
\end{prop}

\proof Consider the 
 basechange map $S_\CC\to S,$
given by $R\to \CC.$ Then, by definition, $\V_1*\V_2|_{S_\CC}$ is 
the parabolic cohomology of a variation of lisse sheaves 
in the sense of \cite{dw03}, Def.~3.1 (namely, the variation 
given by the local system
$\pr_1^*\V_1\otimes \d^*\V_2$). The claim follows then from 
\cite{dw03}, Thm.~3.2.
\Endproof

\subsection{The relation to the 
middle convolution of perverse sheaves.}\label{secperverse}
Let $\ell$ be 
a prime number and let $k$ be an algebraically closed field
of characteristic~$0.$   
In 
$D^b_c(\AA^1_k,{\bar{\QQ}_\ell}),$ define the $*$-convolution (resp. 
the $!$-convolution) as in  the introduction  of this article
or as in \cite{Katz96},
Chap.~2.5. If $\K\in \Perv(\AA^1_k,\bar{\QQ}_\ell)$ has   
the property $\P$ mentioned in the introduction of this article, define the 
middle convolution of perverse sheaves 
$\K\ast_\middle \L \in \Perv(\AA^1_k,\bar{\QQ}_\ell)$ to be 
 $$\K*_\middle\L =\K*\L=\im\left(\K\ast_! \L \to \K\ast_* \L\right).$$
Recall also that for any lisse sheaf $\V$ 
on the punctured affine line $\AA^1_k\setminus \uo,$
one obtains a perverse sheaf $j_*\V[1]$ on $\AA^1_k$ by placing the 
constructible sheaf $j_*\V$ in degree $-1,$ where $j$ denotes the inclusion
of $\AA^1_k\setminus \uo$ into $\AA^1_k.$ 

\begin{prop}\label{verschwindibus2}  
Let $\V_i \in \LS_{\bar{\QQ}_\ell}^\et(U_i),$
 let $\iota_i:U_i \to \AA^1,\, i=1,2,$ denote the 
inclusion maps, and let 
$\K_i:=(\iota_i)_*(\V_1)[1]\in 
\Perv(\AA^1,\bar{\QQ}_\ell).$ Suppose that $\V_1$ is irreducible
and nonconstant.  
Then $\K_1$ has the property $\P$ and 
$$\V_1\ast \V_2= \left(\K_1\ast_\middle\K_2\right)[-1]|_{S},$$
where $S=\AA^1_y\setminus \uo\ast \vo.$
 \end{prop}

\proof That $\K_1$ has the property $\P$
follows from the fact that $\V_1$  is irreducible and not translation
invariant, see \cite{Katz96}, Cor.~2.6.17. By \cite{Katz96},  Prop.~2.9.2,
one has  
\begin{equation}\label{equnten}\K_1\ast_\middle\K_2=R(\overline{\pr}_2)_*\left(j_{!*}(\pr_1^*\K_1\otimes 
\d^*\K_{2})\right),\end{equation} where $j_{!*}$ denotes the middle extension functor
(cf. \cite{Katz96}, Sections~2.3.3 and 2.7, or \cite{KiehlWeissauer}, Chap.~III.5)
associated to 
$j:\AA^1_x\times \AA^1_y \hookrightarrow \PP^1_x\times \AA^1_y$ 
and 
$$ \opr_2:\PP^1_x\times \AA^1_y\to \AA^1_y, \quad \d:\AA^1_x\times \AA^1_y\to \AA^1_{y-x},\, (x,y)\mapsto y-x.$$
 Let $U\subseteq \AA^1_x\times \AA^1_y$ be as in 
Section~\ref{secvondefetale} and let $\tilde{j}:U\to \PP^1\times S$
the inclusion map. Since the divisor $D=\PP^1\times S \setminus U$ is 
smooth over $S$ via $\pr_2,$ 
 the middle extension
$j_{!*}(\pr_1^*\K_1\otimes 
\d^*\K_{2})$ restricted to $\PP^1\times S$ coincides 
with $\tilde{j}_{*}(\pr_1^*\V_1\otimes 
\d^*\V_2)[2]$ (this follows from \cite{Katz96}, Lemma~4.3.8 and 
Cor.~2.8.5.2a). 
It follows therefore from \eqref{equnten} that 
\begin{equation*}\K_1\ast_\middle\K_2|_S=R(\overline{\pr}_2)_*\left(\tilde{j}_*(\V_1\boxtimes \V_2)[2]\right).
\end{equation*} But it follows from Prop.~\ref{propverschwindibus}
that the cohomology sheaves 
$$R^i(\overline{\pr}_2)_*(\tilde{j}_*(\V_1\boxtimes \V_2))
=\H^{i-2}\left(\K_1\ast_\middle\K_2|_S\right)$$
vanish for $i\not=1,$ proving the claim. \Endproof 

\begin{rem} For any algebraically closed field, let $\T_\ell$ denote the category of constructible $\bQl$-sheaves
$\F$ on $\AA^1$ which satisfy the following conditions: 
\begin{enumerate}\item There exist a dense open
subset $j:U\to \AA^1$ such that  $\F|_U$ is irreducible, contained in $\LS_\bQl^\tame (U),$ and  $\F\simeq j_*j^*\F.$ 
\item There are at least two distinct points of $\AA^1$ at which $\F$ 
is not lisse.
\end{enumerate} In \cite{Katz96}, Chap.~5, Katz defines 
a middle convolution functor
$$ \MC_\chi: \T_\ell\To \T_\ell,\quad  \F\mapsto (\F[1]*_\middle 
j_*\L_\chi[1])[-1],$$ where  $\L_\chi$ is the Kummer sheaf associated to 
$\chi:\pi_1^\tame(\GG_m)\to \bQl^\times.$ 
Let $U$ be a dense open subset 
of $\AA^1$ such that $\V=\F|_U$ is lisse.   Let $k$ be a subfield of $\CC.$
It follows then from Prop.~\ref{verschwindibus2} that 
$$ \MC_\chi(\V)=\MC_\chi(\F)|_U.$$ 
By Thm.~\ref{remlab}, the monodromy tuple of $(\MC_\chi(\F)|_U)^\an$ is the 
tuple $\MC_\lambda(\T_{V}).$
\end{rem}

\subsection{Independence of $\ell.$}\label{secfirstpropetalecase}
Let $E$ be a number field and let $\Sigma(E)$ denote the set of 
finite primes of $E.$
Let $S$ be a smooth variety over $R,$ where $R\subseteq \bar{\QQ}$ is  
normal and finitely generated over $\ZZ.$ Suppose that 
for any $\lambda\in \Sigma(E)$ one has given a lisse sheaf 
  $\V_{\lambda} \in \LS_{E_\lambda}^\et(S)$ such that  the following holds:
\begin{itemize}
\item[] Let
 $s:\Spec(\FF_{q})\to S$ be a closed point of $S,$ inducing a map 
$$G_{\FF_{q}}:=\Gal(\bar{\FF}_q/\FF_q)=\widehat{\langle \Frob_s \rangle}\to \pi_1(S,\bar{s}),$$
where $\bar{s}$ is a geometric point which extends $s,$ and  let 
$$ \rho_{\V_\lambda}: \pi_1(S,\bar{s})\To \GL((\V_{\lambda})_{\bar{s}})$$ be the 
monodromy representations of the $\V_\lambda,$ where 
$\lambda\in \Sigma(E).$ Then, for all $\lambda$ of residue characteristic
prime to $q,$ 
 the characteristic 
polynomials of $\rho_{\V_\lambda}(\Frob_s)$ 
have coefficients in $E$ and are
independent of $\lambda.$ 
\end{itemize}

 If this condition is fulfilled, then we say 
that the collection $(\V_\lambda)_{\lambda \in \Sigma(E)}$ is {\it  
$E$-rational and independent of $\lambda$}. The following result is 
close to \cite{Katz96}, Thm.~5.5.4:

\begin{thm}\label{indep}  Suppose that 
for any finite prime $\lambda\in \Sigma(E)$ one has given two nonconstant
irreducible 
pure lisse sheaves
  $$\V_{i,\lambda} \in \LS_{E_\lambda}^\et(U_i),\quad 
 U_i=\AA^1_R\setminus \uo_i\quad (i=1,2),$$
such that the collections $(\V_{i,\lambda})_{\lambda \in 
\Sigma(E)},\,i=1,2,$ are  $E$-rational and 
independent of $\lambda.$ Then the collection
of middle convolutions
$$\left(\V_{1,\lambda}*\V_{2,\lambda}\in \LS_{E_\lambda}^\et(S)\right)_{\lambda \in \Sigma(E)}$$ 
is again $E$-rational 
and independent of $\lambda.$ 
\end{thm}

\proof It suffices
to show that the collection of naive convolutions
$\V_\lambda:=\V_{1,\lambda}*_\naive \V_{2,\lambda}$ is independent of $\lambda$ 
(since the middle convolution $\V_{1,\lambda}* \V_{2,\lambda}$ is the highest-weight-quotient 
of $V_\lambda$ by Prop.~\ref{deligweil}). 
 Let 
$$\F_{\bar{s},\lambda}:=\left(\pr_1^*(\V_{1,\lambda}) \otimes 
d^*(\V_{2,\lambda})\right)|U_{\bar{s}}$$ be the restriction
of the local system $\pr_1^*(\V_{1,\lambda}) \otimes 
d^*(\V_{2,\lambda})$ on $U/S$ (where $U$ is as in 
Section~\ref{secvondefetale}) to the fibre $U_{\bar{s}}$ 
over $\bar{s}\in S(\bar{\FF}_q).$ By the Lefschetz trace
formula, if $\chara(\lambda)\not= \chara(\FF_q),$ then the following holds
(compare to \cite{SGA41/2}, Rapport):
\begin{eqnarray}\label{Lefschetz} 
\sum_{i=0,1,2} (-1)^i \cdot \Trace\left(\Frob_s,H^i_c(X_{\bar{s}},\F_{\bar{s},\lambda})\right)&=&
\sum_{x \in U_s(\FF_{q})} \Trace(\Frob_x,\F_{\bar{x},\lambda}).
\end{eqnarray}
Since $U$ is affine, the 
cohomology group $H^0_c(X_{\bar{s}},\F_{\bar{s},\lambda})$ vanishes
(compare to \cite{FreitagKiehl}, Thm.~9.1). The cohomology group 
$H^2_c(X_{\bar{s}},\F_{\bar{s},\lambda})$ vanishes by Prop.~
\ref{propverschwindibus}. 
Moreover, by the multiplicativity of traces under tensoring, the right
hand side of Equation \eqref{Lefschetz} can be seen to satisfy
\begin{eqnarray*}
\sum_{x \in U_s(\FF_{q})} \Trace(\Frob_x,\F_{\bar{x},\lambda})&=&
\sum_{x \in U_s(\FF_{q})} \Trace\left(\Frob_x,(\V_{1,\lambda})_{\bar{x}}
\right)
\cdot \Trace\left(\Frob_{s-x},(\V_{2,\lambda})_{\overline{s-x}}\,\right).
\end{eqnarray*}
By combining the latter arguments with  Equation \eqref{Lefschetz}, one 
concludes that 
\begin{eqnarray*}\label{Lefschetz3}
\Trace\left(\Frob_s,H^1_c(X_{\bar{s}},\F_{\bar{s},\lambda})\right)&=&
-\sum_{x \in U_s(\FF_{q})} \Trace(\Frob_x,(\V_{1,\lambda})_{\bar{x}})
\cdot \Trace(\Frob_{s-x},(\V_{2,\lambda})_{\overline{s-x}}\,).
\end{eqnarray*}
By the $E$-rationality and independence of 
$\lambda$ of the collections $(\V_{i,\lambda})_{\lambda \in \Sigma(E)},$
this expression is contained in $E$ and is independent of $\lambda.$ 
By definition, it is the trace of $\Frob_s$ on the stalk of the 
naive convolution 
$(\V_{1,\lambda} *_{\naive} \V_{2,\lambda})_{\bar{s}}.$ 
By a similar reasoning,
one sees that  the traces of the powers of $\Frob_s^k$ satisfy
$$
\Trace\left(\Frob_s^k,(\V_{1,\lambda}*_\naive \V_{1,\lambda})_{\bar{s}}\right)
\quad =\quad \quad\quad \quad\quad \quad \quad{} $$
$${}\quad \quad \quad  \sum_{x \in U_s(\FF_{q^{ k}})} \Trace(\Frob_x,(\V_{1,\lambda})_{\bar{x}})
\cdot \Trace(\Frob_{s-x},(\V_{2,\lambda})_{\overline{s-x}}\,)$$
and are thus contained in $E$ and are independent of $\lambda.$ 
Since one can express the coefficients of the characteristic polynomial of 
$\Frob_s$ as linear combinations of the traces of $\Frob_s^k$ (by Newton's
identities), it follows that the collections 
$(\V_{1,\lambda}*_\naive \V_{1,\lambda})_{\lambda \in \Sigma(E)}$ 
and hence $(\V_{1,\lambda}* \V_{1,\lambda})_{\lambda \in \Sigma(E)}$ 
are in fact $E$-rational
and independent of $\lambda.$ 
\Endproof

\begin{thm}\label{thmdet}   Suppose that $R=\ZZ[\frac{1}{N}],$ where 
$N\in \NN_{>0}.$  Let $E$ be a number field and suppose that 
for any finite prime $\lambda\in \Sigma(E)$ one has given two nonconstant
irreducible 
local systems 
  $$\V_{i,\lambda} \in \LS_{E_\lambda}^\et(U_i),\quad {\rm where}\quad 
 U_i=\AA^1_R\setminus \uo_i\quad (i=1,2),$$
such that the collections $(\V_{i,\lambda})_{\lambda \in 
\Sigma(E)},\,i=1,2,$ are  $E$-rational and 
independent of $\lambda.$ Let $S:=\AA^1_R\setminus \uo_1\ast \uo_2$ 
and let 
$V_\lambda:=\V_{1,\lambda}\ast \V_{2,\lambda}\in \LS^\et_{E_\lambda}(S).$
 For $s\in S(\QQ),$ let 
$$ \left(\rho^s_\lambda:G_\QQ:=\Gal(\bar{\QQ}/\QQ)\simeq \pi_1(s,\bar{s})\,\,\To\,\, \Aut(\,(\V_\lambda)_{\bar{s}})\,
\right)_{\lambda\in \Sigma(E)}$$ 
be the induced collection of Galois representations.
Then there exists a finite character 
$\eps:G_\QQ \to E^\times$ and an integer $m\in \ZZ$ such that 
$$\det(\rho^s_{\lambda}) =  \eps \otimes \chi_\ell^m\,,$$
where $\ell=\chara(\lambda)$ and $\chi_\ell$ denotes the
$\ell$-adic cyclotomic character.
\end{thm}

\proof  Let $s=\frac{A}{B}\in S(\QQ)$ with coprime
natural numbers  $A,B\in \ZZ$ and let $\bar{s}$ be a complex point 
extending $s.$ 
The representation $\rho^s_{\lambda}$ 
is unramified outside the prime divisors of $B\cdot N\cdot \ell$
and the finite set of primes $p$ such that $s$ reduces to an element 
in the reduction of $\uo_1\ast \uo_2$ modulo $p.$ 
Thus, by  Thm.~\ref{indep}, the system 
of Galois representations $(\det(\rho^s_{\lambda}))_{\lambda\in \Sigma(E)}$ 
is
a compatible system of $E$-rational Galois representations of $G_\QQ$
in the sense of Serre \cite{Se}.  
  By \cite{Schappacher}, Prop.~1.4 in Chap.~1, 
any such compatible system
arises from an algebraic 
 Hecke character (the Prop.~1.4 in loc.~cit. is a consequence
of Henniart's result on the algebraicity 
of one-dimensional compatible systems, see \cite{Henniart81} and 
\cite{Khare03}). Any Hecke character of $\QQ$ (or, more generally,
a totally real 
field) with values in $E$ is equal to a power of the norm character 
times a finite 
order character with values in $E,$  see \cite{Schappacher}, Chap.~0.3. Therefore,
the associated Galois representation is of the desired form. 
\Endproof

\subsection{Motivic interpretation of the middle
convolution.}\label{secmotiv}

Let $R\subseteq \CC$ be 
normal and  finitely generated over $\ZZ.$
For $i=1,2,$ let  $ U_i=\AA^1_R\setminus \uo_i,$
and let $S:=\AA^1_R\setminus \uo_1\ast \uo_2,$ in the notation of 
Section~\ref{secvondefetale}. 
Suppose that we are given smooth morphisms $\pi_i:X_i\to U_i$ 
and finite groups $G_i$ which are equipped with a 
homomorphism $G_i\to \Aut(X_i/U_i)$ (not necessarily  injective). 
 Let $E$ 
be a number field and 
$P_i\leq E[G_i]$ idempotent elements. 
The homomorphism $G_i\to \Aut(X_i/U_i)$ 
induces a map $E[G_i]\to \End(R^k\pi_{i,!}(E_\lambda))$ for any 
$k\in \NN.$ Thus the idempotent $P_i$ cuts out a subsheaf 
$P_i(R^k\pi_{i,!}(E_\lambda))$ of $R^k\pi_{i,!}(E_\lambda).$
The product $G_1\times G_2$ acts on the scheme
$$\pi: X_1 \boxtimes X_2:=(X_1\times_{U_1} U) \times_U (X_2\times_{U_2} U)
\To U$$
and induces a homomorphism $G_1\times G_2\to \Aut(X_1\boxtimes X_2/U).$
Thus the product idempotent $$P_1\times P_2 \in E[G_1\times G_2]\simeq 
E[G_1]\times E[G_2]$$ 
cuts out subsheaves 
 of the higher direct images $R^k \pi_!(E_\lambda),\,k \in \NN,$
and on the higher direct images 
$R^l(\pr_2\circ \pi)_!(E_\lambda),\, l\in \NN.$ 

Call a constructible $E_\lambda$-sheaf on $U$ 
{\it motivic,} if it is a subfactor of a higher direct image $R^if_!(E_\lambda),$
where $f:X\to U$ is some smooth map.

\begin{thm}\label{thmpropkat}    Let
 $\F_i\in \LS_{E_\lambda}^\et(U_i),\, i=1,2,$ 
be lisse sheaves which are  mixed of weights $\leq k_i$
and which are  of  the following form:
$$ \F_i=P_i\left(R^{k_i}\pi_{i,!}(E_\lambda)\right),
\quad i=1,2,$$
where the $P_i\in E[G_i]$ are idempotents which act 
on $X_i\stackrel{\pi_i}{\to} U_i$ as above.
Assume that $P_i\left(R^{l_i}\pi_{i,!}(E_\lambda)\right)=0$ 
for all $l_i\not= k_i,$ and assume  that 
$\F_1$ or $\F_2$ is irreducible and nonconstant. 
Then 
$$ \F_1\ast_\naive \F_2=P_1\times P_2
\left(R^{k_1+k_2+1}(\pr_2\circ \pi)_{!}(E_\lambda)\right)$$
and  $P_1\times P_2
\left(R^{l}(\pr_2\circ \pi)_{!}(E_\lambda)\right)=0$ for $l\not=k_1+k_2+1.$ 
Moreover, 
$$ W^{k_1}(\F_1)\ast W^{k_2}(\F_2)=W^{k_1+k_2+1}\left(\F_1\ast_\naive \F_2\right).$$
Especially, the naive convolution 
$ \F_1\ast_\naive \F_2$ and the
middle convolution $ W^{k_1}(\F_1)\ast W^{k_2}(\F_2)$ are motivic.  
\end{thm}

\proof The K\"unneth-Formula implies that 
$$ \sum_{i+j=n}R^i\pi_{1,!}(E_\lambda)\otimes 
R^j\pi_{2,!}(E_\lambda)\simeq R^n\pi_!(E_\lambda).$$ 
Thus, 
the assumption 
$P_i\left(R^{l_i}\pi_{i,!}(E_\lambda)\right)=0$ 
for $l_i\not= k_i$ 
implies  that 
\begin{equation} \label{equ28}
P_1[R^{i}\pi_{1,!}(E_\lambda)]\boxtimes 
P_2[R^{j}\pi_{2,!}(E_\lambda)]\simeq 
\left\{\begin{array}{ll}
0 & \textrm{ if } (i,j)\not= (k_1,k_2)\\
P_1\times P_2\left(R^{k_1+k_2}\pi_!(E_\lambda
)\right)&\textrm{ if } (i,j)= (k_1,k_2)\end{array}\right\}
\end{equation}
The Leray spectral sequence 
$E_2=\left(R^p \pr_{2,!}R^q \pi_! (E_\lambda)
\Rightarrow R^{p+q}(\pr_2\circ \pi)_!(E_\lambda)\right)$ 
is a first quadrant 
spectral sequence whose only nonzero entries are at $p=1,2$ (since 
$\pr_2$ is affine, the sheaf $R^{0}\pr_{2,!}$ vanishes). Thus it 
degenerates at $E_2.$ It follows that 
$ R^{n}(\pr_2\circ \pi)_!(E_\lambda)\simeq E^n $ has a descending 
filtration 
$$ R^{n}(\pr_2\circ \pi)_!(E_\lambda)\simeq E_1^n \supseteq E_2^n 
\supseteq 0,$$
where 
$$ E_1^n/E_2^n\simeq R^1\pr_{2,!}R^{n-1}\pi_!(E_\lambda)\quad {\rm and}\quad 
 E_2^n=R^2\pr_{2,!}R^{n-2}\pi_!(E_\lambda).$$ 
It follows therefore from  
\eqref{equ28} and Prop. \ref{propverschwindibus}
that  $ P_1\times P_2\left[R^{n}(\pr_2\circ \pi)_!(E_\lambda)\right]=0$
for $n\not=k_1+k_2+1.$  It follows that 
\begin{eqnarray*} 
\F_1\ast_\naive \F_2&\simeq 
&R^1\pr_{2,!}(P_1(R^{i}\pi_{1,!}(E_\lambda))\boxtimes 
P_2(R^{j}\pi_{2,!}(E_\lambda)))\\
&\simeq & R^1\pr_{2,!}(P_1\times P_2(R^{k_1+k_2}\pi_!(E_\lambda
)))\\
&\simeq &P_1\times P_2\left( R^{k_1+k_2+1}(\pr_2\circ \pi)_!(E_\lambda)\right).
\end{eqnarray*}
By \cite{DeligneWeil2}, the sheaf $R^{k_1+k_2}\pi_!(E_\lambda)$
is  mixed of weights $\leq k_1+k_2.$ 
The exact sequence
$$ 0\to K\to P_1\times P_2\left(R^{k_1+k_2}\pi_!(E_\lambda)\right)\to $$
$$\quad \quad \quad 
 P_1\times P_2\left(W^{k_1+k_2}
(R^{k_1+k_2}\pi_!(E_\lambda))\right)\simeq W^{k_1}(\F_1)\otimes 
W^{k_2}(\F_2)\to 0$$ 
and Prop. \ref{propverschwindibus} imply an exact sequence 
$$ 0\to R^1\pr_{2,!}(K)\to R^1\pr_{2,!}\left(P_1\times P_2\left(R^{k_1+k_2}\pi_!(E_\lambda)\right)\right)$$
$$
\quad\quad\quad\to R^1\pr_{2,!}\left(P_1\times P_2\left(W^{k_1+k_2}
(R^{k_1+k_2}\pi_!(E_\lambda))\right)\right)\to 0.$$ 
Since, again by \cite{DeligneWeil2}, the sheaf 
$R^1\pr_{2,!}(K)$ is mixed of weights $\leq k_1+k_2$ it follows that
$R^1\pr_{2,!}\left(P_1\times P_2\left(W^{k_1+k_2}
(R^{k_1+k_2}\pi_!(E_\lambda))\right)\right)$ is pure of weight $k_1+k_2+1.$ 
It follows then from Prop. \ref{deligweil} (ii) that 
\begin{eqnarray*}
W^{k_1}(\F_1)\ast W^{k_2}(\F_2)&\simeq & W^{k_1+k_2+1}\left(R^1\pr_{2,!}\left(
P_1\times P_2 (R^{k_1+k_2}\pi_!(E_\lambda))\right)\right)\\
&\simeq &W^{k_1+k_2+1}\left(P_1\times P_2(R^{k_1+k_2+1}(\pr_2\circ 
\pi)_!(E_\lambda))\right)\\
&\simeq&W^{k_1+k_2+1}(\F_1\ast_\naive \F_2).
\end{eqnarray*}
\Endproof

\section{Applications to the inverse Galois
problem}\label{secappl}

\subsection{Galois covers and  fundamental groups.}\label{seccovv}

 Let $R$ be a subfield of $\CC,$  let $X$ be 
 a smooth and geometrically irreducible 
variety over $R,$ and let $x$ be a geometric point of X. Any
finite  \'etale Galois cover $f:Y\to X$ with Galois group
$G=G(Y/X)$ corresponds to  a surjective homomorphism   $\Pi_f:
\pi_1^\et(X,x)\to G.$ 
Let $F$ be as in Section~\ref{secetaledef}, let
$\chi:G \rightarrow \GL(V)$ ($V\simeq F^n$)  be a representation,
and let
$\L_{(f,\chi)}\in \LS_F^\et(X)$ 
denote the lisse sheaf associated to the composition
$\chi\circ \Pi_f:\pi_1^\et(X,x)\to \GL(V).$ 
If $X=\AA^1_R\setminus\uo,$ where 
$\uo$ is defined by the vanishing 
of $\prod_{i=1}^r (x-x_i)\in R[x],$ 
let 
$g_{f,\chi}\in \GL(V)^{r+1}$ denote the monodromy tuple which is 
associated to the analytification $\L_{(f,\chi)}^\an$ of $\L_{(f,\chi)}.$
For $r\in \NN_{>0},$ we set
$$ \OO_r(\QQ):=\{P=\{x_1,\ldots,x_r\}\subseteq\bar{\QQ}\mid \,|P|=r
\textrm{ and } P \textrm{ is fixed setwise by 
 } G_\QQ\}.$$

\begin{prop} \label{propreal1} Let $m\in \NN_{>1}$ and let 
$E_\lambda$ be the completion of a number 
field at a finite prime $\lambda$  which admits an embedding 
 $\chi:\ZZ/m\ZZ \hookrightarrow
E_\lambda^\times.$ Let
$$D_m=\ZZ/m\ZZ \rtimes \ZZ/2\ZZ =\langle \sigma,\rho\mid
\rho^m=\sigma^2=1,\,\, \rho^\sigma=\rho^{-1}\rangle$$ denote the
dihedral group of order $2m.$ 
Then the following holds:

 \begin{enumerate}
\item Let $\vo:=\{x_1,x_2\},\, x_1\not= x_2,$ 
where  $x_1,x_2\in \PP^1(\QQ).$
 Then there exists an \'etale Galois cover
$f: {\frak F} \to \PP^1_\QQ\setminus \vo$ with Galois group  isomorphic to
$\ZZ/2\ZZ.$

\item Let $m\in \NN_{>2}$ and let   $\varphi$ denote Euler's $\varphi$-function.
Then there exist infinitely many elements
 $\vo \in \OO_{\varphi(m)}({\QQ})$ such that there is an
\'etale Galois cover $f:{\frak F}\to \AA^1_\QQ \setminus \vo$
 with Galois group isomorphic to
 $\ZZ/m\ZZ=\langle \rho \rangle$ whose monodromy tuple is
$$ \g_{f,\chi}=(\rho^{m_1},\ldots,\rho^{m_{\varphi(m)}},1),$$ where
  $m_i \in \ZZ/m\ZZ^\times.$

\item  Let $m\in \NN_{>2},$ let  $\ZZ/m\ZZ=\langle \rho \rangle,$
and let $\ZZ/2\ZZ=\langle \sigma \rangle.$ 
  Then there exist infinitely many
 $$\vo=\{x_1,\ldots,x_{2+\varphi(m)}\}\in \OO_{2+\varphi(m)}({\QQ})$$ 
for which there exists
an
\'etale Galois cover $f:{\frak F}\to \AA^1_\QQ \setminus \vo$
with Galois group isomorphic to $D_{2m}$ with
$$ \g_{f,\tilde{\chi}}
=(g_1,g_2,\rho^{m_1},\ldots,\rho^{m_{\varphi(m)}},1)\, ,$$
where $g_1,g_2$ are not contained in $\langle \rho
\rangle$ and where $m_i\in \ZZ/m\ZZ^\times.$ Moreover, $x_1$ and $x_2$ can be assumed to
be $\QQ$-rational points.

\item  In the situation of {\rm (iii)} with $m$ odd.
For any $r\geq 2+\varphi(m)$
and (if $r> 2+\varphi(m)$)  
for any choice of elements $x_{3+\varphi(m)},\ldots, x_r\in
\AA^1(\QQ)$ such that 
$$\vo':=\vo \cup \{x_{3+\varphi(m)},\ldots, x_r\} \in \OO_r(\QQ)\, ,$$
there exists  an \'etale Galois cover $f:{\frak F}\to
\AA^1_\QQ \setminus \vo'$ with Galois
group isomorphic to $D_{2m}=\langle \delta \rangle D_m$ (where 
$\delta$ is the central involution)  such that
$$ \g_{f,\tilde{\chi}}=(g_1,g_2,\rho^{m_1},\ldots,\rho^{m_{\varphi(m)}},\delta,\ldots,\delta)
\in D_{2m}^{r+1}. $$ 
\end{enumerate}
\end{prop}

\proof See \cite{Voelklein}, Chap.~7, or \cite{MalleMatzat}, Chap.~I.5.1,
 for (i) and 
(ii). For claim (iii), note that 
the dihedral group $D_m$ is a factor group of 
the wreath product $H=\langle \rho\rangle \wr \langle \sigma \rangle,$
see \cite{MalleMatzat}, Prop.~IV.2.3. Let 
 $\vo=\{x_1,\ldots,x_{\varphi(m)}\}$ be as in (ii) and let 
$$ U:=\{(x,y)\in \CC^2\mid y\not= -{x}{x_i}-x_i^2,\, i=1,\ldots, \varphi(m),
\textrm{ and } 
 x^2\not=4y\}\,.$$
The construction given in \cite{MalleMatzat}, proof of Prop.~IV.2.1, shows that there 
exists an \'etale  Galois cover $\tilde{f}:X\to U$ 
with Galois group isomorphic to $H$  which is defined over 
$\QQ.$
 By factoring out
the kernel of the surjection $H\to D_m,$ one sees that 
there exists 
a Galois cover ${f}':X'\to U,$ with Galois group isomorphic to  $D_m.$
Claim (iii) now follows from restricting the cover $f'$ to a suitable
 (punctured)
line in $U$ which is defined over $\QQ.$
 Claim (iv) follows from (i) and
(iii) by taking  fibre products of covers. \Endproof

\begin{rem}\label{remF} For any Galois cover $f:{\frak F}\to \AA^1_\QQ\setminus \uo$
 as in Prop.~\ref{propreal1} there 
exists an $R=\ZZ[1/r],\, r\in \ZZ,$
 and a Galois cover $f_R: {\frak F}_R\to \AA^1_R\setminus \uo_R$
such that $f$ is the basechange of $f_R$ induced by the inclusion
$R \subseteq \QQ.$ 
\end{rem}

\subsection{Special linear groups as Galois groups.}
 Let $H$ be a profinite group.
One says that $H$ {\em occurs regularly
 as Galois group over $\QQ(t)$} if 
there exists a continuous 
  surjection $\kappa:\Gal(\bar{\QQ(t)}/\QQ(t))\to H $ 
such that $\QQ$ is algebraically closed in the fixed field
of the kernel of $\kappa$. 

\begin{thm}\label{thmrealierung1} Let $\ell$ be a prime, let 
$q=\ell^s\,(s\in \NN),$ and let $n\in \NN.$ Assume that 
$$q\equiv 5\mod 8\quad \mbox{and that}\quad 
n> 6+2\varphi(m)\,,\quad \mbox{ where}\quad  
m:=(q-1)/4 \, .$$
Let 
$$E:=\QQ(\zeta_{m}+\zeta_{m}^{-1},\zeta_4)\, ,$$
where $\zeta_i\,\,(i\in \NN)$ denotes a primitive $i$-th root of unity. Let 
$\lambda$ be a prime of 
$E$ with $\char(\lambda)=\ell$ and let $O_\lambda$ be the valuation ring of 
the completion $E_\lambda.$ 
Then  the 
special linear  group $\SL_{2n+1}(O_\lambda)$ occurs regularly as
Galois group over $\QQ(t).$ \end{thm}

Before giving the proof of the theorem let us mention the
following corollary which immediately follows from reduction modulo
$\lambda$ :

\begin{cor} \label{corsldreia}
The special linear group $\SL_{2n+1}(\FF_q)$  
occurs regularly as Galois
group over $\QQ(t)$ if $$q\equiv 5\mod 8\quad \mbox{and}\quad 
n> 6+2\varphi((q-1)/4)\,.$$ \Endproof
\end{cor}

{\bf Proof of the theorem:} 
First we consider the case where $m\geq 3$  and 
$n=2r-4,$ where  $r\geq  2+\varphi({m}):$ Let 
$$f_1:{\frak F}_1\to \AA^1_\QQ\setminus
\uo_1,\quad \uo_1=\{x_1,\ldots,x_r\}\, ,$$ be a $D_{{2m}}$-cover
as in Prop.~\ref{propreal1}
(iv), where we assume that the points $x_1,\, x_2$ are
$\QQ$-rational. 
Let $\chi_1: D_{{2m}} \hookrightarrow
\GL_2(E_\lambda)$ be an orthogonal embedding of $D_{{2m}}.$ Thus
$$\FFF_1:=\L_{(f_1,\chi_1)}\in \LS^\et_{E_\lambda}(\AA^1_\QQ\setminus
\uo_{1})$$
 is a lisse sheaf of rank two (compare to  Section~\ref{seccovv}
for the notation of $\L_{(f_1,\chi_1)}$). The
monodromy tuple is
$$T_{\FFF_1}=(A_1,\ldots,A_{r+1})\in \GL_2(E_\lambda)^{r+1}\, ,$$
where $A_1,A_2$ are reflections, and $A_3,\cdots,A_{r+1}$ are
diagonal matrices  with eigenvalues
$$(\zeta_{{m}}^{m_1},\zeta_{{m}}^{-m_1}),\ldots,
(\zeta_{{m}}^{m_{\varphi({m})}},\zeta_{{m}}^{-m_{\varphi({m})}}),(-1,-1),
\ldots,(-1,-1)\, ,$$ where 
$\zeta_{{m}}^{m_1=1},\ldots,\zeta_{{m}}^{m_{\varphi({m})}}$ are 
the primitive powers of $\zeta_{{m}}$ (compare to our convention in Section
~\ref{secetaledef} for the notion of an monodromy tuple).

Let $f_2:{\frak F}_2 \to \AA^1_\QQ\setminus \uo_2,\, \uo_2:=\{0\},$ be a double cover
as
in Prop.~\ref{propreal1} (i), let $\chi_2:\ZZ/2\ZZ\hookrightarrow
E_\lambda^\times ,$  and let $$\FFF_2:=\L_{(f_2,\chi_2)}\in
\LS^\et_{E_\lambda}(\AA^1_\QQ\setminus \uo_2)\,.$$ Thus $\FFF_2$ is a
Kummer sheaf with $T_{\FFF_2}=(-1,-1).$

The middle convolution $\FFF_1\ast\FFF_2$ is an \'etale local
system on $\AA^1_F\setminus \uo_1$ (since $\uo_1\ast
\uo_2=\uo_1$). By Prop.~\ref{dimension} and Prop. \ref{dimensione}, 
the rank of
$\FFF_1\ast\FFF_2$ is $2r-4$ and
$$T_{\F_1\ast\F_2}=(B_1,\ldots,B_{r+1})\in \GL_{2r-4}^{r+1}\, ,$$
where  (by Lemma~\ref{lemmonodromy1} and Prop.~\ref{dimension}) the matrices  $B_1,B_2$ are
transvections, the matrices $B_3,\ldots,B_{2+\varphi({m})}$ are
biperspectivities with non-trivial eigenvalues
$$(-\zeta_{{m}}^{m_1},-\zeta_{{m}}^{-m_1}),\ldots,
(-\zeta_{{m}}^{m_{\varphi({m})}},-\zeta_{{m}}^{m_{\varphi({m})}})\quad {\rm (respectively),}$$ the
matrices $B_{3+\varphi({m})},\ldots,B_r$ are unipotent
biperspectivities,  and where the last
matrix is equal to $-1$ by  Lemma~\ref{lemmonodromy2}.

Let $$\GGG_1:=\L_{(g,\xi)}|_{\AA^1_\QQ\setminus \uo_1}\in
\LS^\et_{E_\lambda}(\AA^1_\QQ\setminus \uo_1)$$ be the restriction of
the local system $\L_{(g,\xi)},$ where $g:{\frak G}\to \AA^1_\QQ\setminus
\{x_1\}$ is as in Prop.~\ref{propreal1} (i) and  $\xi : \ZZ/2\ZZ
\hookrightarrow E_\lambda^\times$.  Then the monodromy tuple of
$(\FFF_1\ast \FFF_2)\otimes\GGG_1$ is
$$ T_{(\FFF_1\ast \FFF_2)\otimes\GGG_1}=(-B_1,B_2,\ldots,B_r,-B_{r+1}=1)\,.$$

Let $f_3:{\frak F}_3 \to \AA^1_\QQ\setminus \uo_3$ be a $\ZZ/4\ZZ$-cover
as in Prop.
\ref{propreal1} (ii),  where we have chosen $\uo_3$
such that $\uo_1\ast \uo_3$ is generic in the sense of 
Section~\ref{sectconvdef}. Let
 $\chi_3:\ZZ/4\ZZ \hookrightarrow
E_\lambda^\times$  be an embedding and set $\FFF_3:=\L_{(f_3,\chi_3)}\in
\LS_{E_\lambda}^\et(\AA^1_\QQ\setminus \uo_3).$
 Thus
$$T_{\FFF_3}=(i,-i,1),\quad {\rm where }\quad i:= \zeta_4\,.$$

Then the middle convolution
 $$\V:=((\FFF_1\ast \FFF_2)\otimes\GGG_1)\ast \FFF_3$$ 
is a lisse sheaf
 on $S=\AA^1_\QQ\setminus \uo_1\ast \uo_3.$ The rank of
$\V$ is $4r-7=2n+1$ by Prop.
\ref{dimension}.  
Let $T_{\V}=(C_1,\ldots,C_{2r+1}).$  By Lemma~\ref{lemmonodromy1},
\begin{equation}\label{eq12345}
C_1\sim J(-\zeta_4,2)\oplus_{k=3}^{2r-4}J(-\zeta_4,1)\oplus_{2r-4}^{4r-7}J(1,1)\end{equation} (where $J(\zeta,k)$ denotes 
a Jordan Block of length $k$ with eigenvalue $\zeta$),
$C_2$ is a homology of order four, the
elements $C_3,\ldots,C_{2+\varphi({m})}$ are semisimple
biperspectivities with non-trivial eigenvalues
$$(-i\zeta_{{m}}^{m_1},-i\zeta_{{m}}^{-m_1}),\ldots,
(-i\zeta_{{m}}^{m_{\varphi({m})}},-i\zeta_{{m}}^{m_{\varphi({m})}})\, ,$$ the
elements $C_{3+\varphi({m})},\ldots,C_m$ are biperspectivities
with non-trivial eigenvalues $(i,i),$ and the Jordan forms 
of the matrix $C_{i+r}$ is
the Galois conjugate of the matrix $C_i\, (i=1,\ldots,r).$
Especially, one sees that $\langle C_1,\ldots,C_{2r}\rangle \leq
\SL_{4r-7}(E_\lambda)\times \langle \zeta_4 \rangle.$

%Let $\GGG_2$ be as in Prop.~\ref{real1} and let
%$$ T_{\GGG_2\otimes(\FFF))}
%=(D_1,\ldots,D_{2m+1})\,.$$ Thus $\langle
%D_1,\ldots,D_{2m}\rangle\leq\SL_{4r-7}(E_\lambda).$

It follows from Remark \ref{remF}  that there exists an $N\in \NN_{>0}$ such
that $\V$ extends to a lisse sheaf $\tilde{V}$ 
on $S_{R}=\AA^1_R\setminus (\uo_{1,R}\ast 
\uo_{3,R}),$ where $R=\ZZ[\frac{1}{N\ell}].$
Choose a  $\QQ$-rational 
point $s_0$ of $ S_R$ and hence of $S.$
 By Thm.~\ref{thmdet} and the fact that 
$\tilde{\V}$ is  pure of weight $2$ 
(this follows from 
Prop.~\ref{deligweil}),  
$$\det(\rho_{{\V}})=\det(\rho_{{\V}}|_{\pi_1^{\rm geo}(S,{\bar{s}_0})})
\otimes \chi_\ell^{-(4r-7)}\otimes \epsilon\, ,$$
where $\epsilon:G_\QQ\to \langle
\zeta_4\rangle\subseteq E^\times $ and where $G_\QQ$ is embedded in 
$$\pi_1^\et(S,{\bar{s}_0})=\pi_1^{\rm geo}(S,{\bar{s}_0})\rtimes G_\QQ$$
by the choice of $s_0.$
Thus 
$$\im(\det(\rho_{{\V}(1)}))\leq \langle\zeta_4\rangle \, ,$$
where ${\V}(1)$ stands for the
Tate twist of $\V.$
Let  $$\rho:=\rho_{{\V}}(1)\otimes \delta\;: \;\pi_1^{\rm geo}(S,{\bar{s}_0})\rtimes G_\QQ \To 
\GL_{4r-7}(E_\lambda)\, ,$$ where 
$\delta=\det(\rho_{{\V}(1)})^{-1}.$ It follows that 
\begin{equation}\label{eqgeomm}\im({\rho})\leq \SL_{4r-7}(E_\lambda).
\end{equation}

Let $\FF_q$ be the residue field of $E_\lambda$ and let 
$\bar{\rho}:\pi_1^\et(S,{\bar{s}_0})\to \GL_{4r-7}(\FF_q)$
denote the the residual representation of $\rho.$ Let
 ${H}=\im(\rho)$ and $\bar{H}:=\im(\bar{\rho}).$
Let further ${H}^\geom:=\im(\rho^\geom)$ and 
$\bar{H}^\geom:=\im(\bar{\rho}^\geom).$
By Thm.~\ref{thmirrd} and Prop.~\ref{dimension}, 
$\bar{H}^\geom$ is absolutely 
irreducible. Let $V_1\oplus \cdots \oplus V_l$ be a 
$\bar{H}^\geom$-invariant
decomposition of the underlying module $V:=\FF_q^{4r-7}.$ 
By the eigenvalue structure of the above elements $C_1,C_r\in
{H}^\geom $ given in 
Formula \eqref{eq12345}, the reduction modulo $\lambda$ of
 $C_1$ and $C_r$ does not permute any of the spaces $V_1,\ldots,V_l.$ 
Similarly, for $i=2,\ldots,r,$ the reduction modulo $\lambda$ of
the elements $C_i$ and $C_{r+i}$ fix the spaces $V_1,\ldots,V_l,$
if $\dim(V_i)\geq 3.$ It follows thus that $\dim(V_i)\leq 2.$ 
Since the spaces $V_1,\ldots,V_l$ are permuted transitively by the
irreducibility of $\bar{H}^\geom$ and since the dimension of $V$ 
is odd, we have $\dim(V_i)=1.$ But this is in contradiction to the 
Jordan block of length $2$ which occurs in the reduction modulo $\lambda$ of
 $C_1.$ Therefore,  the group $\bar{H}^\geom$ acts
primitively on $V.$ Thus, the existence of the homology
$C_2$ and  the classification of primitive subgroups of $\GL_n(\FF_q)$ 
 which contain a homology of order greater than $2$ (\cite{Wagner78})
imply that  $\SL_{4r-7}(\FF_q)\leq \bar{H}^\geom.$
By \eqref{eqgeomm}, 
$$\bar{H}^\geom=\SL_{4r-7}(\FF_q)\,.$$

Since $\pi_1^\et(S,{\bar{s}_0})$ is compact, the image of $\rho$ can
be assumed to be contained in the group $\GL_{4r-7}(O_\lambda)$ (see 
\cite{Se})
and hence in the group $\SL_{4r-7}(O_\lambda).$
Since $E_\lambda$ is unramified over 
$\QQ_\ell,$ the residual map 
$\SL_{4r-7}(O_\lambda)\to \SL_{4r-7}(\FF_q)$   is 
Frattini (see \cite{Wei}, Cor.~A). Therefore 
$${H}^\geom=\SL_{4r-7}(O_\lambda)=H\, .$$
The fundamental group
$\pi_1^\et(S,{\bar{s}_0})$ is a factor of
$G_{{\QQ}(t)}$ and the image of $\pi_1^{\rm geo}(S,{\bar{s}_0})$ coincides with the image of $G_{\bar{\QQ}(t)}$ in
$\pi_1^\et(S,{\bar{s}_0}).$ 
Thus
the group
$\SL_{4r-7}(O_\lambda)$ occurs regularly as Galois group over
$\QQ(t).$ This proves the claim in the case $m\geq 3$ and $n=2r-4.$\\

The proof for $m\geq 3$  and $n=2r-3$ uses the following sheaves:
Let
$r\geq 4+\varphi({m}).$ Let $$\FFF_1:=\L_{(f_1,\chi_1)}\in \LS^\et_{E_\lambda}(\AA^1_\QQ\setminus
\uo_1)$$ be defined as above, where we assume 
that the points $x_{\varphi(m)+3}$ and 
$x_{\varphi(m)+4}$ coincide with $i,-i$ and 
that $x_1,x_2$ and $x_r$ are rational points. 
 By a suitable tensor operation,
one obtains a local system 
$\FFF_1'\in \LS^\et_{E_\lambda}(\AA^1_\QQ\setminus
\uo_1)$  whose
monodromy tuple is
$$T_{\FFF_1'}=(A_1',\ldots,A_{r+1}')\in \GL_2(E_\lambda)^{r+1}\, ,$$
where $A_1',A_2'$ are reflections, and $A_3',\cdots,A_{r+1}'$ are
diagonal matrices  with eigenvalues
\begin{multline}(\zeta_{{m}}^{m_1=1},\zeta_{{m}}^{-m_1}),\ldots,
(\zeta_{{m}}^{m_{\varphi({m})}},\zeta_{{m}}^{-m_{\varphi({m})}}),\\
\quad (i,i),\,(-i,-i),\,(-1,-1),\ldots,(-1,-1),(1,1))\,.\end{multline}
Let $\FFF_2=\L_{(f_2,\chi_2)}$ be the Kummer 
sheaf as above and let $\FFF_3':=\FFF_2.$ 
Let 
$\GGG'=\L_{(g',\xi')},$ where $g': {\frak G}'\to \AA^1_\QQ\setminus\{x_1,x_r\}$ 
is the double cover ramified at $x_1$ and $x_r$ and 
$\xi'$ is the embedding of the Galois group 
into $E_\lambda.$
Now continue as above, using the sheaf 
$$\V=((\FFF_1'\ast \FFF_2)\otimes \GGG'|_{\AA^1_\QQ\setminus \uo_1})\ast 
\FFF_3'$$ which has rank  $4r-5.$

The case where $m=1$ follows from  the same arguments as above
using the dihedral group $D_3$ instead of $D_m.$
\Endproof

\section{Appendix (by S. Reiter and M. Dettweiler): A new motivic local system of $G_2$-type}

Let $\J(n_1,\ldots,n_k)$ denote a unipotent matrix in Jordan 
canonical form which decomposes into blocks of length $n_1,\ldots,n_k$
and let $\zeta_3\in \bQl$ denote a fixed 
primitive third root of unity. Let us call 
a constructible $\bQl$-sheaf on a scheme $X$  to be 
{\it motivic} is it is a subfactor 
of a higher direct image sheaf of a morphism $Y\to X.$ 
We obtain the following result:\\

\noindent{\bf Theorem:}
{\it  Let $k\subseteq \CC$ be an algebraically closed 
field and 
let  $x_1,x_2,x_3\in \AA^1(k)$ be pairwise distinct.
 Then there exists a motivic 
lisse sheaf $\HHH$  on $\PP^1_k\setminus \{x_1,x_2,x_3,\infty\}$ 
with the following properties:
\begin{enumerate}
\item The Jordan canonical forms of the 
local monodromies of $\HHH$ at $x_1,x_2,x_3,\infty$ are as follows (resp.):
$$ \J(2,2,1,1,1),\,\,\J(2,2,1,1,1),\, \,
\diag(1,\zeta_3,\zeta_3,\zeta_3,\zeta_3^{-1},\zeta_3^{-1},\zeta_3^{-1}),\,\, 
\J(3,3,1)\,.$$ 
\item The monodromy of $\HHH$ is dense in 
$G_2(\bQl).$  
\item The lisse sheaf $\HHH$ is not cohomologically rigid in the sense 
of \cite{Katz96}, Chap.~5.\\
\end{enumerate}}

\proof 
Let  
$X=\PP^1_k\setminus \{x_1,x_2,x_3,\infty\}.$ Remember from
Section~\ref{seccovv} that 
to any \'etale Galois cover $f:Y\to X$ with Galois 
group $G$ and to any homomorphism $\alpha:G\to \GL_n(\bQl)$ we have associated 
a lisse sheaf $\L_{f,\alpha}.$ Let us also fix generators
$\alpha_1,\ldots,\alpha_4$ of the topological 
fundamental group $\pi_1(X(\CC))$
and 
generators $\gamma_1,\gamma_2$ of $\pi_1(\GG_m(\CC))$ as in the
Definition of the monodromy tuple (Definition~\ref{deff}).
Let $\L_i=\L_{f_i,\alpha}\in 
\LS_\bQl^\et(X),\, 
i=1,2,3,$ be as follows:
\begin{itemize}
\item
The Galois cover $f_1:Y_1\to X$ is the cover of $X$
defined by $$y^3=(x-x_1)(x-x_3).$$  
The Galois group of $f_1$ is 
$C_3=\langle \sigma\mid \sigma^3=1 \rangle.$
The  homomorphism $\alpha: C_3\to \bQl^\times$ is given by 
sending $\sigma$ to $\zeta_3.$ 
The monodromy tuple of the analytification $\L_1^\an \in \LS(X(\CC))$ 
is $\g_1=(\zeta_3,1,\zeta_3,\zeta_3).$ 
\item The Galois cover $f_2:Y_2\to X$ is the cover of $X$
defined by $$y^3=(x-x_2)(x-x_3)$$  and $\alpha$ is as above. 
The monodromy tuple of $\L_2^\an\in \LS(X)$ 
is equal to $\g_2=(1,\zeta_3,\zeta_3,\zeta_3).$ 
\item The Galois cover $f_3:Y_3\to X$ is the cover of $X$
defined by $y^3=(x-x_3)$ and $\alpha$ is as above. 
 The monodromy tuple of $\L_2^\an\in \LS(X)$ 
is equal to $\g_3=(1,1,\zeta_3,\zeta_3^{-1}).$
\end{itemize}
Let $\chi:\pi_1(\GG_m)\to \bQl^\times$ be the character, which sends 
the generator $\gamma_1$ of $\pi_1(\GG_m)$ to $\zeta_3$ and let 
$\chi^{-1}$ be its dual character. 
Consider the following sequence
of tensor operations and middle convolutions: 
\begin{equation*} \HHH:=\L_3^{-1}\otimes(\MC_{\chi^{-1}}(\L_3
\otimes(\MC_\chi(\MC_{\chi^{-1}}(\L_1)\otimes \MC_{\chi^{-1}}(\L_2))))) \,
\in 
\LS_\bQl^\et(X),
\end{equation*} where $\L_3^{-1}$ denotes the dual sheaf of $\L_3.$ 
It follows from Prop.~\ref{dimension}~(i) that
then analytification $\HHH^\an \in  
\LS_\bQl(X(\CC))$ is of the form 
\begin{equation*} \HHH^\an=(\L_3^{-1})^\an\otimes(\MC_{\chi^{-1}}(\L_3^\an
\otimes(\MC_\chi(\MC_{\chi^{-1}}(\L_1^\an)\otimes \MC_{\chi^{-1}}(\L_2^\an))))).\end{equation*} 
Using the explicit recipe for the computation of the middle convolution
of local systems
given in Thm.~\ref{remlab}, 
one finds explicit matrices 
for the monodromy tuple 
 $h=(h_1,\ldots,h_4)$ of $\HHH^\an.$ These are given in Table~\ref{Table1}.
\begin{table}%[ht]
\begin{center}
\begin{small}
$$ h_1=\left(\begin{array}{ccccccc}
     1& -3& \zeta_3 - 1& 0& \zeta_3 - 4 &0& 2\zeta_3 + 4\\
     0& 3\zeta_3 + 1& 2\zeta_3 + 1& 0& 2\zeta_3 + 1& -2\zeta_3 - 1& 0\\
     0& -3\zeta_3 &-2\zeta_3& 0 &-2\zeta_3 - 1& 2\zeta_3 + 1 &0\\
     0& 3\zeta_3 + 3& \zeta_3 + 2 &1& \zeta_3 + 2& -\zeta_3 - 2& 0\\
     0 &3\zeta_3 + 6& 3& 0& 4& -3& 0\\
     0& 3\zeta_3 + 3& \zeta_3 + 2& 0& \zeta_3 + 2& -\zeta_3 - 1& 0\\
     0& 6& -2\zeta_3 + 2 &0 &-2\zeta_3 + 2& 2\zeta_3 - 2& 1
\end{array}\right)$$
$$ h_2=\left(\begin{array}{ccccccc}
   1& 0& 0& 0& 0& 0& 0\\
   \zeta_3 - 1& 1& 0& 2\zeta_3 + 1& 0& 0& 0\\
     3& 0 &1 &-2\zeta_3 - 1& -3& 0& 2\zeta_3 + 4\\
     0& 0 &0 &1& 0& 0& 0\\
     0 &0& 0& 0& 1& 0& 0\\
     0& 0& 0 &0& 0& 1 &0\\
     0 & 0& 0 &0& 0& 0 &1\end{array}\right),$$ $$
 h_3=\left(\begin{array}{ccccccc}
     \zeta_3 &0 &0 &0& 0& 0& 0\\
     0& \zeta_3 &0& 0& 0& 0& 0\\
     0& 0& \zeta_3& 0& 0& 0& 0\\
     \zeta_3 + 2& 0& 0& -\zeta_3 - 1& 0& 0& 0\\
     0& \zeta_3 + 2& 0& 0& -\zeta_3 - 1& 0& 0\\
     0 &3\zeta_3 + 3& \zeta_3 + 2& 0& 0& -\zeta_3 - 1& 0\\
     0& 0& 0& 0& \zeta_3 - 1& 0& 1\end{array}\right)\,.$$
\end{small}

\caption{The monodromy generators of $\HHH$} \label{Table1}
\end{center}
\end{table}
 The claim on the local monodromy of $\HHH$ follows then from the Jordan 
forms of the matrices $h_1,\ldots,h_4.$

Since the monodromy
tuples of $\MC_{\chi^{-1}}(\L_i^\an),\,i=1,2,$ are contained 
in the group $\Sp_2(\bQl)=\SL_2(\bQl)$ and generate 
irreducible subgroups, the monodromy tuples
of the tensor product 
 $\MC_{\chi^{-1}}(\L_1^\an)\otimes \MC_{\chi^{-1}}(\L_2^\an)$
can be seen to generate an irreducible subgroup
of 
the orthogonal group $\SO_4(\bQl).$
An iterative application of Thm.~\ref{thmirrd}  implies then that
the local system $\HHH$ is irreducible. 

Since the elements 
of  the monodromy tuple of 
 $\MC_{\chi^{-1}}(\L_1^\an)\otimes \MC_{\chi^{-1}}(\L_2^\an)$
are contained in $\SO_4(\bQl),$ it follows from
Poincar\'e duality that
the elements of the monodromy tuple of 
$\HHH^\an$  are contained  in the group $\SO_7(\bQl).$ 
By a computation using
MAGMA \cite{magma}, one can check that the matrices 
stabilize a one-dimensional subspace of the third exterior power
$\Lambda^3(\bQl^7)  .$ Thus,  by the results of \cite{Asch},
the image of the monodromy representation
$\rho^\an$ of $\HHH^\an$ 
is contained in $G_2(\bQl).$  By the classification
of bireflection groups given in \cite{GuralnickSaxl}, Thm.~7.1 and Thm.~8.3, 
the Zariski closure of the 
image of $\rho^\an$ can be seen to coincide
with $G_2(\bQl).$ Since $\pi_1^\et(X)$ is the 
profinite closure of $\pi_1(X(\CC)),$ 
the Zariski closure of  
 $\rho^\an$ coincides 
with the Zariski closure of the monodromy representation 
of $\HHH.$ 
 It follows from the numerical criterion 
of physical rigidity given in 
\cite{Katz96}, Thm. 1.1.2, and the structure of the local monodromy of $\HHH^\an$ 
 that
$\HHH^\an$ is not 
physically rigid. Therefore, $\HHH$ not cohomologically rigid in the sense
of \cite{Katz96}, Chap.~5. 

To prove that $\HHH$ is motivic we can work over 
$\ZZ[\zeta_3,\frac{1}{3}]$ instead of $k$ and use an 
iterative application 
of the K\"unneth-formula and the motivic interpretation
of the middle convolution given in Thm.~\ref{thmpropkat}. 
\Endproof

\noindent{\bf Remark:} 
Let $\zeta_6$ denote a primitive 
$6$-th root of unity and let 
$\g_1:=(\zeta_6,1,\zeta_6,\zeta_6),$  $\g_2=(1,\zeta_6,\zeta_6,\zeta_6)$
and   $\g_3=(1,1,\zeta_6,\zeta_6^{-1}).$ Let 
$\L_{\g_i},\, i=1,2,3,$ be the local systems on $X(\CC)$ associated
to $\g_i$ and let 
 $\chi:\pi_1(\GG_m(\CC))\to \bQl^\times$ be the character, which sends 
the generator $\gamma_1$ of $\pi_1(\GG_m)$ to $\zeta_6.$ Then
the sequence
of tensor operations and middle convolutions
\begin{equation*} \L_{\g_3}^{-1}\otimes(\MC_{\chi^{-1}}(\L_{\g_3}
\otimes(\MC_\chi(\MC_{\chi^{-1}}(\L_{\g_1})\otimes \MC_{\chi^{-1}}(\L_{\g_2})))))
\end{equation*} leads to a local system on $X(\CC)$ whose local
monodromy coincides with the one of $\HHH$ but whose monodromy
is Zariski dense in $\SO_6(\bQl).$ 

    \bibliographystyle{plain}
                    \bibliography{$HOME/Biblio/p}

\end{document}